\documentclass[twoside,12pt]{amsart}
\usepackage[active]{srcltx} 
\usepackage[all]{xy}
\CompileMatrices
\usepackage{amsfonts}
\usepackage{amssymb}
\usepackage{amsthm}
\usepackage{amsmath}
\usepackage{ifthen}
\usepackage{url}
 \newtheorem{thm}{Theorem}[section]
 \newtheorem{prop}[thm]{Proposition}
 \newtheorem{cor}[thm]{Corollary}
 \newtheorem{lem}[thm]{Lemma}
  
 \newtheorem{conj}[thm]{Conjecture}
 \newtheorem*{bthm}{Theorem}
 \newtheorem*{bprop}{Proposition}

\theoremstyle{definition}
\newtheorem{defn}[thm]{Definition}

\newtheorem{ex}[thm]{Example}

\theoremstyle{remark}
\newtheorem{rem}[thm]{Remark}

\font\russ=wncyr10  1
\def\sha{\hbox{\russ\char88}}

\newcommand{\fp}{\ifmmode {\mathbb{F}_p}\else$\mathbb{F}_p$\ \fi}
\newcommand{\zp}{\ifmmode {\mathbb{Z}_p}\else$\mathbb{Z}_p$\ \fi}
\newcommand{\z}{\mathbb{Z}}
\newcommand{\zpMod}{\ifmmode\mbox{$\zp$-Mod}\else$\zp$-Mod \fi}
\newcommand{\Mod}{\ifmmode\mbox{$\Lambda$-Mod}\else$\Lambda$-Mod \fi}
\renewcommand{\mod}{\ifmmode\mbox{$\Lambda$-mod}\else$\Lambda$-mod
\fi}
\newcommand{\La}{\ifmmode\Lambda\else$\Lambda$\fi}

\newcommand{\Hom}{{\mathrm{Hom}}}

\newcommand{\Tor}{{\mathrm{Tor}}}

\newcommand{\pd}{{\mathrm{pd}}}
\newcommand{\cd}{{\mathrm{cd}}}
\newcommand{\rk}{{\mathrm{rk}}}
\newcommand{\gl}[1]{{\mathrm{gl}(#1)}}

\renewcommand{\H}{\mathrm{H}}

\newcommand{\M}{\ifmmode {\frak M}\else${\frak M}$ \fi}
\newcommand{\m}{\ifmmode {\frak m}\else$\frak m$ \fi}
\newcommand{\mh}{\ifmmode {\frak m}(H)\else${\frak m}(H)$ \fi}
\newcommand{\p}{\ifmmode {\frak p}\else${\frak p}$\ \fi}
\renewcommand{\P}{\ifmmode {\frak P}\else${\frak P}$\ \fi}
\newcommand{\e}{\ifmmode {\mathcal{E}}\else$\mathcal{E}$ \fi}

\newcommand{\T}{\mathcal{ T}}
\renewcommand{\O}{\mathcal{ O}}
\newcommand{\G}{\ifmmode {\mathcal{G}}\else${\mathcal{G}}$\ \fi}
\renewcommand{\d}{\ifmmode {\mathcal{ D}}\else${\mathcal{D}}$\ \fi}
\newcommand{\A}{\ifmmode {\mathcal{A}}\else${\mathcal{ A}}$\ \fi}

\newcommand{\F}{\mbox{$\mathcal{F}$}}
\renewcommand{\projlim}[1] {{\lim\limits_{\stackrel{\displaystyle
\longleftarrow}{#1}}}}

\newcommand{\dirlim}[1]
{{\lim\limits_{\stackrel{\displaystyle \longrightarrow}{#1}}}}
\renewcommand{\in}{\ \epsilon\ }

\newcommand{\kl}{[\![}
\newcommand{\kr}{]\!]}
\newcommand{\Qp}{\ifmmode {{\Bbb Q}_p}\else${\Bbb Q}_p$\ \fi}
\newcommand{\qp}{\ifmmode {{\Bbb Q}_p}\else${\Bbb Q}_p$\ \fi}
\newcommand{\Q}{\ifmmode {\Bbb Q}\else${\Bbb Q}$\ \fi}

\newcommand{\Ind}{\mathrm{Ind}}
\newcommand{\coker}{\mathrm{coker}}

\newcommand{\gr}{\mathrm{gr}}
\newcommand{\Sel}{\mathrm{Sel}}

\def\sectionnam{\@empty}
\def\subsectionnam{\@empty}
\newcommand{\Section}[1]{\section{#1}}

\newcommand{\Subsection}[1]{\subsection{#1}}

\begin{document}

\title{Characteristic elements in noncommutative Iwasawa theory\\\mbox{ }\\}%
\author{Otmar Venjakob}%
\address{Universit\"{a}t Heidelberg\\ Mathematisches Institut\\
Im Neuenheimer Feld 288\\ 69120 Heidelberg, Germany.} \email{otmar@mathi.uni-heidelberg.de}
\urladdr{http://www.mathi.uni-heidelberg.de/\textasciitilde otmar/}



\date{\today}%

 \maketitle

%
%
%
%
%
%


%
\PUSH{intro.tex}%

\thispagestyle{empty}

Let $p$ be a prime number, which, for simplicity, we shall always assume odd. In the Iwasawa theory of an elliptic curve $E$ over a number field $k$ one has to distinguish between curves which do or do not admit complex multiplication (CM). For  CM-elliptic curves  their deep arithmetic properties and the link between their Selmer group and  special values of their Hasse-Weil $L$-functions are not only described by the (one-variable) main conjecture corresponding to the cyclotomic $\zp$-extension $k_{cyc}$ of $k,$ but also by the (two-variable) main conjecture corresponding to the extension $k_\infty=k(E_{p^\infty})$ which arises by adjoining the $p$-power division points $E_{p^\infty}$ of $E.$ Moreover, both conjectures are proven by Rubin \cite{rubin} in the case that $k$ is imaginary quadratic and $E$ has CM by the ring of integers $\mathcal{O}_k$ of $k.$\\
 Also for non-CM elliptic curves one would like to at least formulate   a main conjecture over the trivialzing extension $k_\infty,$ but for lack of both an algebraic as well as analytic $p$-adic $L$-function this has not been achieved.
The aim of this paper is to establish, under certain conditions, the existence of an {\em algebraic} $p$-adic $L$-function, viz as an element of the first $K$-group $K_1(\La_\T)\cong\La_\mathcal{T}^\times/[\La_\mathcal{T}^\times,\La_\mathcal{T}^\times]$ of a localization $\La_\T$ of the usual Iwasawa algebra $ \La=\La(G)$ of the Galois group $G=G(k_\infty/k).$ Here, for a ring $R,$ we denote by $R^\times$ its group of units.  By the Weil-pairing, $k_{cyc}$ is contained in $k_\infty$ and we put $H=G(k_\infty/k_{cyc})$ and $\Gamma=G(k_{cyc}/k).$ Furthermore we write $\mh$ for
the kernel of the canonical surjective ring homomorphism
  \[\psi_H:\La(G) \twoheadrightarrow
\mathbb{F}_p\kl \Gamma\kr\] and we set    $\mathcal{T}:=\La\setminus\mh.$ Since $\mh$ is a (completely) prime ideal $\T$ is obviously a multiplicative subset of $\La=\La(G).$ But in contrast to commutative ring theory it is not clear at all whether one can localize a ring at a given multiplicative set.  Thus, from a technical point of view the following theorem is the starting point of this paper

\begin{bthm}[Theorem \ref{ore-thm}] Let  $G$ be isomorphic   the
direct product $H\times \Gamma$ with $H$  a uniform group (i.e.\ a pro-$p$-group without element of order $p$ and such that $[H,H]\subseteq H^p$). Then the
multiplicative set $\mathcal{T}$ satisfies the left and right Ore condition
for $\La(G).$  In particular, the localization $\La_\mathcal{T}$ exists and is a Noetherian regular local ring with global dimension $\gl {\La_\mathcal{T}}<\gl \La.$
\end{bthm}

The proof relies on filtered ring techniques and makes decisive use of Lazard's fundamental work \cite{la} which we review in section \ref{lazard}. 
The   significance of the existence of $\La_\mathcal{T},$ which we will assume always henceforth,  results from the observation that  the relative $K$-group $K_0(\La,\La_\mathcal{T})$ is closely linked with a nice category of $\La$-torsion modules. By $\mbox{\rm $\La(G)$-mod}^H$  we denote the full subcategory of \mbox{\rm $\La(G)$-mod} consisting of those $\La(G)$-modules which are finitely generated over the subalgebra $\La(H)$ of $\La=\La(G).$ We will see below that such modules play an important role in our arithmetic applications.  

\begin{bprop}[Proposition \ref{t-torsion}]
There are a canonical isomorphisms of groups \[K_0(\La(G)\mbox{-mod}^H)\cong K_0(\Lambda(G),\La(G)_{\T})\cong \La_\mathcal{T}^\times/[\La_\mathcal{T}^\times,\La_\mathcal{T}^\times]\La^\times
.\]
\end{bprop}

This identification enables us to define characteristic elements $F_M $ in $\La_\mathcal{T}^\times,$ unique only modulo $[\La_\mathcal{T}^\times,\La_\mathcal{T}^\times]\La^\times,$ for modules $M$ in $ \La(G)\mbox{-mod}^H$ requiring  that   it corresponds via the above isomorphism to the class $[M]\in K_0(\La(G)\mbox{-mod}^H).$  
From the definition of $\mathcal{T}$ it is clear that the projection $\pi_H:\La(G)\to \La(\Gamma)$ extends to a ring homomorphism $\pi_H:\La(G)_\mathcal{T}\to Q(\Gamma),$ where $Q(\Gamma)$ denotes the field of fractions of $\La(\Gamma).$ Thus we obtain a commutative ``descent" diagram of $K$-groups with exact rows
\[  \xymatrix{
  K_1(\La) \ar[r]\ar@{->}^{(\pi_H)_\ast}[d] & K_1(\La_\mathcal{T})\ar[r]\ar@{->}^{(\pi_H)_\ast}[d] & K_0(\La,\La_\mathcal{T}) \ar[r]\ar@{->}^{(\pi_H)_\ast}[d]& 0 \\
  K_1(\La(\Gamma)) \ar[r]\ar@{=}[d] & K_1(Q(\Gamma))\ar[r]\ar@{=}[d] & K_0(\La(\Gamma),Q(\Gamma)) \ar[r]\ar@{=}[d]& 0 \\
  {\La(\Gamma)^\times \ar[r] }& Q(\Gamma)^\times\ar[r] & Q(\Gamma)^\times/\La(\Gamma)^\times\ar[r] &
  1.
}\]

which can be used to define the evaluation of $F_M$ at `` $0$ " or using twisting (see section \ref{twisting}) more generally at certain $p$-adic representations, see section \ref{evaluation} for details. Thus let $\rho: G\to GL(V)$ be a continuous linear representation on a finite dimensional vector space $V$ of dimension $m$ over a finite extension $K$ of $\qp$ with ring of integers $\O.$ We choose a $G$-invariant $\O$-lattice $T\subseteq V$ and define a generalized $G$-Euler characteristic as follows \[\chi (G,V,M):=\prod_i
(\#\mbox{Tor}_i^{\La(G)}(T,M))^{(-1)^{i}}\] provided that all Tor-groups are finite. One checks easily that this is independent of the choice of $T.$\\
Since the class of $M$ in $K_0(\Lambda(G),\La(G)_{\T})$ can also be described using virtual objects (see section \ref{swan}) which should be considered as universal Euler characteristic of $M$ and which behave well under change of rings one obtains immediately the following relation between the characteristic element of a module $M$ and its $G$-Euler characteristic, which is very important for descent arguments in arithmetic applications, e.g.\ if one wished to link the desired main conjecture with the Birch and Swinnerton-Dyer conjecture.

\begin{bthm}[Theorem \ref{EC-evaluation}] Assume that $M$ belongs to $\mbox{\rm $\La(G)$-mod}^H.$ 
Then, if \linebreak $\chi(G,V,M)$ is finite,  $F_M(\rho)$ is defined 
and we have
\[\chi(G,V,M)=|F_M(\rho)|_p^{-[K:\qp]}\] for any choice of $F_M.$ In particular, $F_M(\rho)$ is non-zero.
\end{bthm}

Having settled these purely algebraic properties we want to verify its usefulness in the study of Selmer groups over $p$-adic Lie extensions. Beneath $k_\infty=k(E_{p^\infty})$ we also consider the    false-Tate-curve-case, i.e.\ $k_\infty$ arises as trivializing extension of a $p$-adic representation which is analogous  to that of the local Galois representation  associated with a Tate elliptic curve. More precisely, we assume that $k$ contains the group  $\mu_p$ of $p^\mathrm{th}$ roots of unity  and then $k_\infty$  is obtained by adjoining to $k_{cyc}$  the $p$-power roots of an element in $k^\times$ which is not a root of unity. By Kummer theory, its Galois group is isomorphic to the semidirect product $G(k_\infty/k)\cong\zp(1)\rtimes \zp$ where the action  is given by the cyclotomic character.\\
We write  \[X(k_{cyc})=X_f(E/k_{cyc})=(\Sel_{p^\infty}(E/k_{cyc}))^\vee\] for the the Pontryagin dual of the Selmer group $\Sel_{p^\infty}(E/k_{cyc})$ of $E$ over   $k_{cyc}$ and similarly $X(k_\infty)$ over $k_\infty.$ Then $X(k_{cyc})$ is a finitely generated $\La(\Gamma)$-module and in case it is torsion (as it is always conjecturally)  we denote by $f_{X(k_{cyc})}$ its (classical) characteristic polynomial. \\ 
The strongest confirmation that our algebraic $p$-adic $L$-function bears indeed arithmetic information relies on the fact that the characteristic element $F_{X(k_\infty)}$ of the dual of the Selmer group over $k_\infty$ specializes,  up to some local Euler factors,  to the characteristic polynomial $f_{X(k_{cyc})}$ of the dual of the Selmer group over $k_{cyc}.$

\begin{bthm}[Theorems \ref{tate-case}, \ref{GL2-case}] 
Assume that in the false-Tate-curve- or $GL_2$-case $X(k_{cyc})$ is a torsion $\La(\Gamma)$-module with
vanishing $\mu$-invariant. Then $X(k_\infty)$ is finitely
generated over $\La(H)$ and we have  \[\pi_H(F_{X(k_\infty)})\equiv
f_{X(k_{cyc})} \cdot \prod_{\nu\in
\mathfrak{M}} \mathcal{P}_\nu(E(p)/k)\;\;\mathrm{  mod}\; \Lambda(\Gamma)^\times, \] where $\mathcal{P}_\nu(E(p)/k)$ are  local Euler factors, see  section \ref{local-euler}, while $\mathfrak{M}$ denotes a certain finite set of primes of $k$ which ramify in the extension $k_\infty/k,$ see subsections \ref{tate-curve} and \ref{GL2}.
\end{bthm}

Beneath the techniques developed above this result relies heavily on the vanishing of higher $H$-homology groups of $X(k_\infty)$ and   descent calculations done by J. Coates, P. Schneider and  R. Sujatha \cite{cs2} in the $GL_2$-case and by Y. Hachimori and the author \cite{hachi-ven} in the false Tate curve case.
By evaluation of the above formula at the trivial representation we reobtain under the above conditions the   well known determination of the $G$-Euler characteristic which can now be interpreted as the value at `` $0$ " of the characteristic element of the dual of the Selmer group.

\begin{bthm}[Corollary \ref{BSD}]
In the situation of the theorem and assuming that the
$G$-Euler characteristic $\chi(G,X(k_\infty))$ of
$X(k_\infty)$ is finite let $F_{X(k_\infty)}\in \La_\mathcal{T}$ be a
characteristic element of $X(k_\infty).$ Then $F_{X(k_\infty)}(0)$ is
defined (and non-zero) and it holds that
\[\chi(G,X(k_\infty))=|F_{X(k_\infty)}(0)|_p^{-1}=\rho_p(E/k)\times\prod_{v\in
\mathfrak{M}} |L_v(E,1)|_p.\]
Here, $\mathfrak{M}$ denotes the same set of primes as above; for the definition of the $p$-Birch-Swinnerton-Dyer constant $\rho_p(E/k)$ see subsection \ref{tate-curve}.
\end{bthm}


These results  encourage us to suggest a possible shape of a main conjecture in section \ref{mainconj} involving our definition of characteristic element for  the dual of the Selmer group. For a precise conjecture one needs a good guess for the $\epsilon$-factors, Deligne period, etc.\ and since there is very little empirical material available at the moment this  seems a very delicate and subtle point. But a similar recipe as in \cite{coates91} for the conjectural $p$-adic $L$-functions  of motives over $\mathbb{Q}$ with respect to the cyclotomic $\zp$-extension should generalize to our situation.

\vspace{0,3cm}
\textsc{Acknowledgements.}
\smallskip

This   is a condensed version of an article which has been submitted as Habilitationsschrift to the Fakult{\"a}t f{\"u}r Mathematik und Informatik der Ruprecht-Karls-Universit{\"a}t Heidelberg.
I would like to thank  heartily Kay Wingberg for his confidence and advice, which accompanied my research during the last years. To John Coates I express my warmest  gratitude for his great interest and encouragement which pushed this work  permanently. Also I am indebted to Sujatha   for many valuable comments on an earlier version of this article. I am very grateful to Peter Schneider for his invitation to participate in the SFB programme "Geometrische Strukturen" in M{\"u}nster for one week; our intensive joint discussions at that time helped me some months later to discover the definition of the Ore set $\T$ which is the right one for proving the existence of characteristic elements. Hearty thanks go to Yoshitaka Hachimori for stimulating conversations during his stay in Heidelberg. David Burns is heartily acknowledged for illustrating the equivariant Tamagawa number conjecture and the concept of virtual objects to me. Also I thank Denis Vogel for allowing me to discuss his examples of non-principal reflexive ideals in this paper.\\

 \vspace{0,3cm}
\textsc{General Notation and Conventions}
\smallskip
\begin{enumerate}
\item In this paper, a ring $R$ is always associative and with a unit element. When we are talking about properties related to $R$ like being ``Noetherian", an ``ideal", a
``unit", of a certain ``global dimension", etc.\ we always mean
the left {\em and} right property if not otherwise stated. But by
an $R$-module we usually mean left $R$-module  (not a bi-module).
\item By $R^\times$ we denote the group of units, i.e.\ of right and left invertible elements.
\item By $\pd_R(M)$ we denote the
projective dimension of a $R$-module $M$ while $\mathrm{gl}(\Lambda)$ denotes the global
dimension of $\Lambda.$
\item By a {\em local} ring $R$ we mean a ring in which the non-units
form a proper ideal, which is then automatically maximal as left,
right and two-sided ideal. Equivalently, $R$ has both a unique
left and a unique right maximal ideal, which amounts to the same
as the quotient $R/J(R)$ of $R$ by its Jacobson radical being a
skewfield.
\item A filtration $F_\bullet
R:=\{F_nR|n\in\mathbb{Z}\}$ of a ring we shall  always assume is
indexed by $\mathbb{Z},$  increasing, 
exhaustive and
separated. We write $\gr R=\bigoplus_{n\in\mathbb{Z}} F_nR/F_{n-1}R$
for its associated graded ring. A similar notation and convention
is used for filtered left $R$-modules.  
\item   For a  discrete  (resp.\ compact) $\zp$-module $N$
with continuous action by some profinite group $G$,
$N^\vee=\Hom_{\zp,cont}(N,\Q_p/\zp)$ is the compact (resp.\
discrete) Pontryagin dual of $N$ with its natural $G$-action. If
$N$ is  $p$-divisible, $T_p(N)=\Hom(\Q_p/\zp,N)= \projlim{i}\;
{_{p^i}}N$ denotes the Tate module of  $N,$ where ${_{p^i}}N$ denotes the kernel of the multiplication by  $p^i.$ For $G=G_k$ the  absolute Galois group of a number field or a local field  $k,$ we define the $r$th Tate twist   of $N$ by
 $N(r):=N\otimes_\zp T_p(\mu)^{\otimes
 r}$ for $r\in\mathbb{N}$ and $N(r):=N\otimes_\zp \Hom(T_p(\mu)^{\otimes
 r},\zp)$ for $-r\in \mathbb{N},$ where $\mu$ denotes the $G_k$-module of all
 roots of unity and by convention
 $T_p(\mu)^{\otimes 0}=\zp$ with trivial $G$-action. Finally,
we set $N^\ast:=\varinjlim_i\Hom({_{p^i}}N,
 \mu_{p^\infty})=T_p(N)^\vee(1).$
\end{enumerate}

\POP

\PUSH{structure.tex}%
\Section{$p$-adic Lie groups}\label{Lie}

Since the towers of number fields studied in non-commutative Iwasawa theory form a $p$-adic Lie extension in this subsection we recall basic facts about  $p$-adic Lie groups. The reader who is familiar with this topic may skip this section or only glance at it for notational reasons.

There is a famous characterization of $p$-adic analytic groups
due to Lazard \cite{la} (see also \cite{dsms2} 9.36):
{\em A topological group $G$  is a compact $p$-adic Lie group if
and only if
 $G$ contains a normal open uniformly powerful
pro-$p$-subgroup of finite index.}
Let us briefly recall the definitions:  A pro-$p$-group $G$ is
called {\em powerful}, if  $[G,G]\subseteq G^{p}$ for odd $p,$
respectively $[G,G]\subseteq G^{4}$ for $p=2,$ holds. Here for any prime $p$ and natural number $n$ we write $G^{p^n}$ for the subgroup of $G$ which is generated by all elements of the form $g^{p^n},$ $g\in G.$ A
(topologically) finitely generated powerful pro-$p$-group $G$ is
{\em uniform} if it has no
element of order $p$ (see \cite{dsms2}, p. 62). For instance, for $p\geq n+2,$ the group $Gl_n(\zp)$
has no elements of order $p,$ in particular, $GL_2(\zp)$ contains
no elements of finite $p$-power order if $p\geq 5$ (see
\cite{howson} 4.7). It follows that
  all the congruence kernels of $GL_n(\zp),SL_n(\zp)$ or
$PGL_n(\zp)$ are uniform pro-$p$-groups for $p\neq 2.$ 

%

We should mention that  Lazard himself did not use the notation of powerful or uniform groups. Instead, he formulated the above characterization of $p$-adic analytic groups using the notation of {\em $p$-valuable} groups, i.e.\  complete $p$-valued groups of finite rank, see  \cite[thm. III.3.1.7]{la}.  
 We recall that a group $G$ is called {\em $p$-valued} \cite[III, Def.\ 2.1.2]{la} if it possesses a $p$-valuation, i.e.\  a function $\omega: G\to (0,\infty]$ satisfying the following axioms for all $g$ and $h$ in $G:$ \begin{enumerate} 
\item $\omega(1)=\infty,$ and $1/(p-1)<\omega(g)<\infty$ for $g\neq 1,$
\item $\omega(gh^{-1})\geq \min\{\omega(g),\omega(h)\},$ 
\item $\omega(g^{-1}h^{-1}gh)\geq \omega(g)+\omega(h)$ and 
 \item $\omega(g^p)=\omega(g)+1.$
\end{enumerate}

In particular, it follows that $G_\nu:=\{g\in G|\omega(g)\geq
-\nu\}$ and $G_{\nu^+}:=\{g\in G|\omega(g)> -\nu\}$ are normal
subgroups for each $\nu$ in $\mathbb{R}.$ A $p$-valued group
$(G,\omega)$ is said to be {\em complete} if $G=\projlim\nu G/G_\nu.$
Putting \[\gr\: G:=\bigoplus_{\nu\in\mathbb{R}} G_\nu/G_{\nu^+}\] we
obtain a graded Lie algebra (the Lie bracket being  induced by the
group commutator)  over the graded ring $\fp[\pi_0]=\gr\:\zp$ with
$\pi_0$ in degree $-1$ (the action of $\pi_0$ on $\gr\: G$
corresponds to taking the $p^\mathrm{th}$ power of an element of $G,$ see
\cite[pp.\ 464-465]{la}). In fact, $\gr\: G$ is
 free as a $\gr\:\zp$-module and its rank is called the {\em rank} of $(G,\omega)$ (\cite[III
 2.1.3]{la}). In particular, a  $p$-valued group has no element of
 order $p.$ The class of 
 $p$-valuable groups  is closed under taking closed subgroups and forming finite
 products.
 We should remark that by \cite[III 3.1.11]{la} a $p$-valuable group always admits a
 $p$-valuation which  has rational
 values,  more precisely in a discrete subset of $\Q.$ Henceforth
 we assume that $\omega$ is of the latter sort.
If a compact $p$-adic group $G$ is $p$-valued it is automatically
pro-$p$ (\cite[III 3.1.7]{la}.  
As mentioned above, Lazard's original characterization now reads as:
{\em A topological group is a compact
$p$-adic Lie group (of finite dimension $d$) if and only if it
contains a normal open  $p$-valuable subgroup of finite rank
$d$ (\cite[III 3.1.3/7/9, 3.4.5]{la}). In this case the dimension
 $\dim G$ (of the underlying $p$-adic manifold) and the rank of $G$ coincide (loc.\ cit.).}

As explained in the remark after \cite[lem.\ 4.3]{sch-teit2} the relation between  $p$-valuable and uniform pro-$p$-groups (for $p\neq 2$ for simplicity) is as follows (see also the Notes at the end of chapter 4 in \cite{dsms2}): A $p$-valuable group $G$ is called {\em $p$-saturated} \cite[III.2.1.6]{la} if $G$ has a $p$-valuation $\omega$ with the property that any $g\in G$ with $\omega(g)>p/(p-1)$ is a $p^\mathrm{th}$ power. If moreover, $G$ has a minimal system of (topological) generators $g_1,\ldots,g_d$ such that $\omega(g_i)+\omega(g_j)> p/(p-1)$ for any $1\leq i\neq j\leq d,$ we say that $G$ is {\em strongly $p$-saturated} as this property implies that every commutator in $G$ is a $p^\mathrm{th}$ power. In particular,  every strongly $p$-saturated $p$-valued group $G$ is powerful, and thus uniform since it does not have an element of order $p.$ Conversely, uniform groups even allow a $p$-valuation $\omega$ with $\omega(g_i)=1$ for every $g_i$ of an arbitrary minimal system of generators $g_1,\ldots, g_d.$ Since any element in $G_2$ is a $p^\mathrm{th}$ power \cite[le.\ 4.10]{dsms2} one has:
{\em For $p\neq 2,$ a $p$-valuable group $G$ is uniform if and only if it is strongly $p$-saturated.}

\Section{The Iwasawa algebra - a review of Lazard's work}\label{iwasawa-algebra}

{\em Throughout this section we make the   following assumption:
Let $\O$ be a commutative Noetherian local ring which is complete in its
$\m$-adic topology, where $\m$ is the maximal ideal. We assume
that $\kappa=\O/\m$ is a finite field of characteristic $p,$ in
particular $\O$ is compact.}

 In our applications $\O$ is usually
the ring of integers in a finite extension field of $\qp$ or a
finite field.

We denote by $\Lambda=\Lambda(G)$ the {\em Iwasawa algebra}  of a
compact $p$-adic Lie group $G,$ i.e.\ the completed group algebra
of $G$ over $\O$
 $$\Lambda(G)=\O\kl G\kr =\projlim{U}\O[G/U],$$
where $U$ runs through the open normal subgroups of $G.$ For a
good treatment of basic properties of $\La,$ some of which we
recall below, we refer the reader to \cite[V{\S}2]{nsw}. First we
should mention that 
{\em $\La$ is a semi-local ring, it is local if and
only if $G$ is pro-$p.$}

The global dimension $\mathrm{gl}\La(G)$ equals
$\cd_pG+\mathrm{gl}\O$ where $\cd_p$ denotes $p$-cohomological
dimension \cite{brumer}. By a result of Serre \cite{serre65} $\cd_pG$ is finite if and only if
$G$ does not contain an element of order $p.$

The whole deeper structure theory of $\La$ relies on the following
observation which is essentially due to Lazard \cite{la}.

\begin{thm}{\rm (Lazard)}\label{lazard}
Assume that $G$ is a $p$-valuable (hence compact $p$-adic Lie) group and, in addition to our general assumptions, that $\O$
is  a finite field or a discrete valuation ring (DVR).  Then $\O \kl
G\kr$ possesses a complete separated and exhaustive (increasing)
filtration $F_\bullet \La$ such that
\begin{eqnarray*}
\gr\: \O \kl G\kr &\cong & \left\{\begin{array}{rl}
  \kappa[X_0,X_1,\ldots,X_d]&  \mbox{ if $\O$ is a DVR} \\
  \kappa[X_1,\ldots,X_d] & \mbox{ if $\O=\kappa$ is a finite field}
\end{array}\right.
\end{eqnarray*}
are isomorphisms of graded rings. Here $d=\dim G,$ and the grading
on the polynomial ring is given by assigning to each variable a
certain strictly negative integer degree.
\end{thm}

For $\O=\zp$ this reformulation of Lazard's result is \cite[prop.\
7.2]{css}. Before we extend it to a bigger class of rings $\O$ we
restate Lazard's original results more precisely. First note that
any valuation $\omega$ of $G$ extends to a filtration also called $\omega$ of the Iwasawa
algebra $\zp\kl G\kr$ and with respect to this filtration
 Lazard \cite[Ch.\ III 2.3.3/4]{la} has established
a canonical isomorphism
\[\gr_\omega\:\zp\kl G\kr\cong U(\gr_\omega\: G).\]  Using this results he shows
that  one can  describe the elements of $\zp\kl G\kr$ as certain
power series in non-commuting variables. For this purpose we fix an {\em ordered basis} of
$G,$ i.e.\ a sequence of elements $g_1,\ldots , g_d\in
G\setminus\{1\}$ such that the elements $g_iG_{\omega{g_i}^+}$
form a basis of $\gr\: G$ as an $F_p[\pi_0]$-module.
Then any $\lambda\in\zp\kl G\kr$ has a unique convergent expansion
\[\lambda=\sum_{\alpha\in \mathbb{N}^d} \lambda_\alpha (g_1-1)^{\alpha_1}\cdot\ldots\cdot (g_d-1)^{\alpha_d},\]
with $\lambda_\alpha\in \zp$ for all $\alpha$ and conversely all
such series converge in $\zp\kl G\kr$ (\cite[III 2.3.8]{la}, see
also \cite[{\S}4]{sch-teit2}). In the language of Lazard this means
that $ \La$ is complete-free (\cite[I 2.1.17]{la}) with topological
basis
\[ \mathbf{b}_\alpha:=(g_1-1)^{\alpha_1}\cdot\ldots\cdot
(g_d-1)^{\alpha_d}, \alpha\in \mathbb{N}^d.\] In other words
there is an isomorphism of filtered (in particular topological)
$\zp$-modules
\[\zp\kl G\kr \cong \prod_\alpha\zp \mathbf{b}_\alpha,\] where the
filtration of the left hand side is the product filtration  with
 $\zp \mathbf{b}_\alpha$ being isomorphic to the  filtered $\zp$-module
 $\zp (r_\alpha)$ which is isomorphic to $\zp$ but with  shifted filtration by $r_\alpha=\sum_{i=1}^d \omega(g_i).$

%
%
%
%

Now observe that there is a canonical isomorphism
of topological $\O$-algebras  \[\O\kl G\kr\cong \O\otimes_\zp \zp\kl
G\kr.\] Thus the tensor product filtration on the right hand side
induces a filtration on $\O\kl G\kr$ such that there is an isomorphism of filtered (thus in
particular topological) $\O$-modules \[\O\kl G\kr\cong
\prod_\alpha\O \mathbf{b}_\alpha\] and an isomorphism of graded
$\gr\:\O$-algebras
\[\gr\:\O\kl G \kr\cong U(\gr\:\O\otimes_{\gr\:\zp}\gr\:
G),\] where $U(\gr\:\O\otimes_{\gr\:\zp}\gr\: G)$ denotes the enveloping
algebra of the Lie algebra $ \gr\:\O\otimes_{\gr\:\zp}\gr\: G.$
Furthermore, it is always possible to replace the valuation
$\omega$ by a valuation $\omega'$ with values in $\Q$ such that
$(G,\omega')$ is also $p$-valued and such that $\gr_{\omega'} G$
is an abelian Lie algebra. In this case, we obtain an isomorphism
\[\gr_{\omega'}\O\kl G \kr\cong (\gr\:\O)[X_1,\ldots,X_d]\] where $d$ is the
dimension of $G.$ Now the theorem follows immediately  in the general case.
Also we obtain the following Proposition where by
  $\m_G$ we denote the maximal ideal of $\O\kl G\kr.$

\begin{prop}\label{uniform}
Let $G$ be a uniform group and assume  that $\O$ is as in the above theorem.   Then,  with respect to the $\m_G$-adic filtration there is an
isomorphism of (possibly non-commutative) graded $\gr\:\O$-algebras
\[\gr_{\m_G}\:\O\kl G \kr\cong U(\gr\: G\otimes_{\gr\:\zp}\gr\:\O).\]
\end{prop}

\begin{rem}
Lazard uses filtrations indexed by the positive real numbers in
general. But in \cite[proof of 7.2/3]{css} it is explained that in
our situation one always can rescale the filtration on $\zp\kl
G\kr$ to get one indexed by $\z.$ Of course, this extends to general
$\O\kl G\kr.$
\end{rem}

Now, the strategy is to use techniques from the theory of filtered
rings and modules to fully exploit  Lazard's Theorem and to
derive further properties of $\La$ and of its module category.
Most of these methods can be found in the book \cite{li-o}. 

\begin{cor}
Assume that $\O$ is  finitely generated as a $\zp$-module and let
$G$ be a compact $p$-adic Lie group. Then $\O\kl G\kr$ is
Noetherian.
\end{cor}

\begin{proof}
Using Lazard's characterization of $p$-adic Lie groups and the
fact that for an open subgroup $H$ of $G$ the module $\La(G)$ is finitely
generated over $\La(H)$ we may assume that $G$ itself is
$p$-valued. Now the result follows since  $\gr\:
\La(G)\cong(\gr\:\O)[X_1,\ldots ,X_d]$ is Noetherian and $\La(G)$ is
complete (\cite[II{\S}1 prop. 3]{li-o}).
\end{proof}

Using similar techniques one  shows that for a $p$-valued compact
$p$-adic Lie group the ring $\O\kl G\kr$ has no zero divisor if
$\O$ is either a finite field or a DVR:

\begin{cor}
In the situation of the Theorem $\O\kl G\kr$ is an integral domain.
\end{cor}

In the mixed characteristic case  A. Neumann \cite{Ne} shows this for all
torsionfree pro-$p$ $p$-adic Lie groups (in the case $\O=\zp$) by
different means, which unfortunately do not generalize e.g.\ to
the case of a finite field $\O= \kappa.$

\begin{prop}{\rm (Neumann)}
Assume that $\O$ is regular with mixed characteristic and let $G$
be a pro-$p$-group without any element of order $p$ and such
$\O\kl G\kr$ is Noetherian (e.g.\ if $G$ is a pro-$p$ $p$-adic Lie
group). Then $\O\kl G\kr$ has no zero divisor.
\end{prop}

The proof is completely analogous to \cite[thm 1]{Ne} using a
theorem of Walker since the finiteness of the global dimension is
well known by Brumer's  result.

\POP

\PUSH{k-theory0.tex}%
\Section{Virtual objects and some requisites from $K$-theory}
\label{swan}

We begin recalling Swan's construction of relative $K$-groups. For any ring $\La$ we denote by $\mathcal{P}(\La)$ the category of finitely generated projective $\La$-modules. For any homomorphism of rings $\phi:\La\to\La',$ the relative $K$-group $K_0(\La,\La')$ is defined by generators and relations, as follows. Consider triples $(M,N,f)$ with $M,N\in \mathcal{P}(\La),$ $f:\La'\otimes_\La M\cong  \La'\otimes N.$ For brevity, let $M'=\La'=\La'\otimes_\La M,$ etc. A morphism \[(\mu, \nu):(M_1,N_1,f_1)\to (M_2,N_2,f_2)\] consists of a pair of maps $\mu\in\Hom_\La(M_1,M_2),\; \nu\in \Hom_\La(N_1,N_2),$ such that \[ \nu'\circ f_1=f_2\circ \mu':M'_1\to N'_2.\] We write $(M_1,N_1,f_1)\cong (M_2,N_2,f_2)$ if both $\mu$ and $\nu$ are isomorphisms. A short exact sequence of triples is a sequence 
\[\xymatrix@1{
   {\ 0 \ar[r] } &  {\ (M_1,N_1,f_1)\ar[r]^{(\mu_1,\nu_1)} } &  {\ (M_2,N_2,f_2) \ar[r]^{(\mu_2,\nu_2)} } &  {\ (M_3,N_3,f_3) \ar[r] } &  {\ 0 } 
}\]
such that each pair $(\mu_i,\nu_i)$ is a morphism, and where the sequences of $\La$-modules
\[\xymatrix@1{
   {\ 0 \ar[r] } &  {\ M_1 \ar[r]^{\mu_1} } &  {\ M_2 \ar[r]^{\mu_2} } &  {\ M_3 \ar[r] } &  {\ 0 } 
}\]
and similarly for $N_i$ with $\nu_i$ are exact. Now $K_0(\La,\La')$ is defined as the free abelian group generated by all isomorphism classes of triples, modulo the relations \[(L,N, gf)=(L,M,f)+(M,N,g)\] and for each short exact sequence as above 
\[(M_2,N_2,f_2)=(M_1,N_1,f_1)+(M_3,N_3,f_3).\] This relative $K$-group fits into the following exact sequence of groups 
\[\xymatrix@1{
   {\ K_1(\La) \ar[r] } &  {\ K_1(\La') \ar[r]^{\delta} } &  {\ K_0(\La,\La') \ar[r]^{\lambda} } &  {\ K_0(\La) \ar[r] } &  {\ K_0(\La'),  } 
}  \]
where the map $\delta$ is defined by $\delta(f)=[\La^n,\La^n,f]$ for $f\in GL_n(\La'),$ while the map $\lambda$ is given by $\lambda([M,N,f])=[M]-[N],$ and where the brackets denote classes of triples in $K_0(\La,\La')$ and $K_0(\La),$ respectively.

Next we consider the special case where $\La'$ arises as localisation $\La_\mathcal{T}$ of a Noetherian regular ring $\La$ without zero-divisors by an Ore set $\T.$
Recall that a multiplicative closed subset $\mathcal{T}$ of a ring $R$ is
said to satisfy the {\em right Ore condition} if, for each $r\in
R$ and $s\in \mathcal{T},$ there exist $r'\in R$ and $s'\in \mathcal{T}$ such that
$rs'=sr'.$ If $R$ is Noetherian, then the right Ore condition
guarantees that the right localisation $R_\mathcal{T}$ of $R$ at $\mathcal{T}$ exists.
There is an analogous left version of this and we say that $\mathcal{T}$ is
an {\em Ore set} if it satisfies both the left and right Ore
condition. In this case the left and right localisation are
canonically isomorphic and thus identified and called localisation
of $R$ at $\mathcal{T}.$ A good reference for (classical) localisation is the book
\cite[Ch.\ 2]{mc-rob}.

We say that a $\La$-module $M$ is $\T$-torsion if $\La_\T\otimes_\La M=0$ and we denote by $\La\mbox{-mod}_{\T-\mathrm{tor}}$ the full subcategory of $\La\mbox{-mod}$ consisting of all $\T$-torsion modules. Now the category $\La_\T\mbox{-mod}$ can be identified with the quotient category of $\La\mbox{-mod}$ with respect to the Serre subcategory  $\La\mbox{-mod}_{\T-\mathrm{tor}}.$ Using that by regularity of the rings $\La$ and $ \La_\T$ their $G$- and $K$-theory coincide the localisation exact sequence of  $K$-theory looks like
\[\xymatrix@1{
   {\ K_1(\La) \ar[r] } &  {\ K_1(\La') \ar[r]^-{\delta} } &  {\ K_0(\La\mbox{-mod}_{\T-\mathrm{tor}}) \ar[r]^-{\lambda} } &  {\ K_0(\La) \ar[r] } &  {\ K_0(\La'),  } 
}  \]
where the map $\lambda$ is induced by the inclusion of categories while the map $\delta$ is defined by $\delta(f)=[\coker(f)]$ for $f\in GL_n(\La')\cap M_n(\La)$ and noting that any element of $ GL_n(\La')$ is a product of the form $fg^{-1}$ with $f,g\in GL_n(\La')\cap M_n(\La),$ see \cite{bass}.  

It is well known and  follows from the $5$-lemma that there is a canonical isomorphism of groups 
\begin{eqnarray}
 K_0(\La\mbox{-mod}_{\T-\mathrm{tor}}) &\cong &K_0(\La,\La_\T)
\end{eqnarray}
once we have established a map commuting with the $\delta$'s and $\lambda$'s. Since in our applications both $\La$ and $\La_\T$ are regular rings we only describe it in this situation, for simplicity. First note that under this assumptions the maps $\lambda$ are both trivial because $[\La]\in \z[\La]\cong K_0(\La)$ is mapped to $0\neq[\La_T]\in \z[\La_\T]\cong K_0(\La_\T).$ Now let $M$ be in $\La\mbox{-mod}_{\T-\mathrm{tor}}$ and choose a (finite) projective, thus free resolution 
\[\xymatrix@1{
   {\ F^\bullet=F^\bullet(M):} & { \cdots \ar[r] } & {  F^i\ar[r] } &{ \cdots \ar[r] } & {\ F^2 \ar[r] } &  {\ F^1 \ar[r] } &  {\ F^0 \ar[r] } &  {\ 0 } 
}\]
of $M.$ After tensoring with $\La_\T$ this becomes an acyclic resolution $F^\bullet_\T$ of free $\La_\T$-modules by assumption on $M.$ By $F^+=F^+(M)$ and $F^-=F^-(M)$ we denote the even and odd sum \[\bigoplus_{i\;\mathrm{ even}} F^i \;\;\mbox{ and } \;\; \bigoplus_{i\,\mathrm{ odd}} F^i, \] respectively, and similarly for $F^\bullet_\T.$ Now we choose successively sections in order to obtain  a map $\phi:F^+_\T\to F^-_\T$ and we define the image of $M$ in $K_0(\La,\La_\T)$ as \[M\mapsto [F^+,F^-,\phi].\] We leave it to the reader to check that this is independent of the choice of the resolution $F^\bullet$ as well as of the map $\phi,$ for more details in a slightly different context see also \cite[lem.\ 1.1.3,1.2.2,1.2.3, thm.\ 1.2.1]{burns}. One also has to check that our map is additive on short exact sequences in order to see that the above assignment really induces a map on classes $[M]\in K_0(\La\mbox{-mod}_{\T-\mathrm{tor}}),$ cf.\ (thm.\ 1.2.7, loc.\ cit.). For an example of this construction, see \ref{ex-char} (i).

Note that under the assumption that both $\La$ and $\La_\T$ are local, i.e.\ $K_1(\La_\T)\cong \La_\T^\times/[\La_\T^\times,\La_\T^\times]$ and similarly for $\La$ by \cite[ex.\ 1.6]{srinivas}, the group $K_0(\La,\La_T)$ is generated by triples of the form $[\La,\La,f]$ with $f\in \La_T^\times\cap\La.$ One sees immediately that $[\La,\La,f]\mapsto [\coker(f)]$ induces the inverse map $ K_0(\La,\La_\T)\cong K_0(\La\mbox{-mod}_{\T-\mathrm{tor}}). $

Now we shall
describe the image of $M$ in $K_0(\La,\La_\mathcal{T})$ using the
concept of virtual objects due to Deligne \cite{d} and further
developed by Burns and Flach \cite{bf} and also used by Huber and Kings \cite{hu-ki} in their approach to non-commutative Iwasawa theory.
Let $R$ be a associative ring  with unit. 
The category of {\em virtual objects} $V(R)$ is
a Picard category, i.e.\ a groupoid ( a category in which all
morphisms are isomorphisms) equipped with a bifunctor
$(L,M)\mapsto L\boxtimes M$ satisfying certain conditions, see
\cite[{\S}2.1]{bf} for details, in particular it has a unit object $
\mathbf{1}_{V(R)}$ unique up to unique isomorphism. Furthermore
$V(R)$ comes equipped with a functor
\[ [-]: (D^p(R),is)\to V(R),\] where $D^p(R)$ denotes the  category
of perfect complexes (as full triangulated subcategory of the
derived category $D^b(R)$ of the homotopy category of bounded
complexes of $R$-modules) and $(D^p(R),is)$ denotes the
subcategory of isomorphisms. For the construction of that functor
and a list of its compatibility properties we refer the reader to
(prop. 2.1, loc. cit.), here we only mention that $[-]$ commutes
with the functors $R'\otimes_R-:D^p(R)\to D^p(R')$ and
$R'\otimes_R-:V(R)\to V(R')$ induced by any ring extension $R\to
R'$ (prop. 2.1 d), loc.\ cit.) and that for each exact triangle in
$D^p(R)$
\[\xymatrix@1{
   {\ \Sigma=\Sigma(u,v,w):X \ar[r]^-{u} } &  {\ Y \ar[r]^-{v} } &  {\ Z \ar[r]^-{w} } &  {\ X[1]
   }}
    \]
    there is a nonempty set $[\Sigma]$ of isomorphisms
   $\phi:[Y]\cong [X]\boxtimes [Z]$ in $V(R);$ in case $R$ is a
   regular ring $[\Sigma]$ consists of precisely one element,
   $\phi_\Sigma$ say
(prop. 2.1, loc.\ cit.). Finally,  one has for $X\in D^p(R)$ such
that all $\H_i(X)\in D^p(R),$ too, a canonical isomorphism
\begin{eqnarray}\label{det-homology}
[X]\cong
\raisebox{-0,2cm}{$\stackrel{\displaystyle\boxtimes}{\scriptstyle
i\in \z}$} [\H_i(X)]^{(-1)^{i-1}},
\end{eqnarray} see
(loc.\ cit.\ (9)).\\
For a commutative local ring the pair $(V(R),[-])$ is equivalent
to the determinant functor  in the sense of Knudsen and Mumford
\cite{km} taking values in the category of graded line bundles on
$Spec(R)$ (this is the best way to think about virtual objects).
Finally we should mention that the fundamental groups of the
Picard category $V(R)$ ($\pi_0(V(R))$ is the group of isomorphism
classes of objects of $V(R)$ while
$\pi_1(V(R))=\mathrm{Aut}_{V(R)}(\mathbf{1}_{V(R)})$) are
canonically isomorphic to the $K$-groups $K_0(R)$ and $K_1(R)$ 
of $R$ ({\S}2.3, loc.\ cit.). Also the relative $K$-group $K_0(R,R')$
for a ring homomorphism $R\to R'$ can be realized as fundamental
group of a Picard category: Let $\mathcal{P}$ be the Picard
category  with unique object $\mathbf{1}_\mathcal{P}$ and
$\mathrm{Aut}_\mathcal{P}(\mathbf{1}_\mathcal{P})=0.$ Following
\cite[(20)]{bf} we define $V(R,R')$ to be the fibre product
category $V(R)\times_{V(R')} \mathcal{P}.$ Thus objects of
$V(R,R')$ consists of pairs $(M,\lambda)$ with $M\in V(R)$ and
$\lambda: R'\otimes_R M\to\mathbf{1}_{V(R')}$  an isomorphism in
$V(R').$ In analogy with prop. 2.4 (loc.\ cit) we obtain an
isomorphism \[K_0(R,R')\cong \pi_0(V(R,R'))\] where $[M,N,f]$ is mapped to $([M]\boxtimes [N]^{-1},[f]\boxtimes \mathrm{id}_{[R'\otimes_R N]^{-1}}).$

Now we return to our previous situation and let $M$ be in
$\mbox{$\La$-mod}_{\mathcal{T}\mathrm{-tor}}.$ Then the
associated complex $M[0]$ concentrated in degree $0$ belongs to
$D^p(\La)$ and its image in $D^p(\La_\mathcal{T})$ is given by an
acyclic complex. Thus the isomorphism $\La_\mathcal{T}\otimes_\La
M [0]\to 0$ in $D^p(\La_\mathcal{T})$
induces an isomorphism \[\lambda_M:\La_\mathcal{T}\otimes_\La [M
[0]]=[\La_\mathcal{T}\otimes_\La M[0]]\to
[0]=\mathbf{1}_{V(\La_\mathcal{T})}\] and hence the pair
$([M[0]],\lambda_M)$ is an element in $V(\La,\La_\mathcal{T})$. We
denote its class in $K_0(\La,\La_\mathcal{T})$
 by $\mathrm{char}_\La(M)$  and call it the
{\em characteristic class} of $M.$

From the analog of \cite[rem.\ in 2.7]{bf} we obtain the following

\begin{prop}\label{ident}
Under the above identifications $K_0(\La\mbox{-mod}_{\T-\mathrm{tor}})=K_0(\La,\La_\T)=\pi_0(V(\La,\La_\T))$ the different classes associated to $M\in\La\mbox{-mod}_{\T-\mathrm{tor}}$ coincide \[[M]=[F(M)^+,F(M)^-,\phi]=\mathrm{char}_\La(M).\]
\end{prop}

Since $K_1$ of a local ring can be calculated via the Dieudonn\'{e}
determinant and the relative $K$-group $K_0(\La,\La_\mathcal{T})$
 fits into the following short exact sequence
\[ \xymatrix{
  K_1(\La) \ar[r]\ar@{=}[d] & K_1(\La_\mathcal{T})\ar[r]\ar@{=}[d] & K_0(\La,\La_\mathcal{T}) \ar[r]\ar@{=}[d]& 0 \\
  {\La^\times/[\La^\times,\La^\times]\ar[r]} & {\La_\mathcal{T}^\times/[\La_\mathcal{T}^\times,\La_\mathcal{T}^\times]\ar[r]} & {\La_\mathcal{T}^\times/[\La_\mathcal{T}^\times,\La_\mathcal{T}^\times]\La^\times\ar[r] }&
  1,} \]
we can consider $\mathrm{char}_\La(M)$ also as an element $F_M
[\La_\mathcal{T}^\times,\La_\mathcal{T}^\times]\La^\times\in
\La_\mathcal{T}^\times/[\La_\mathcal{T}^\times,\La_\mathcal{T}^\times]\La^\times.$
Then we call any such choice $F_M\in\La_\mathcal{T}^\times$ a {\em characteristic}
element of $M.$

\begin{rem}\label{ex-css2}
Let $G$ be a pro-$p$ $p$-adic Lie group without element of order $p$ and denote by $Q(G)$ its skew field of fractions of $\La(G).$ Then  the above applies   to the ring homomorphism $\La(G)\to Q(G)$ and in this case the relative $K$-group $K_0(\La(G),Q(G))$ describes nothing else than the Grothendieck group  of the full torsion subcategory $\mbox{\rm $\La(G)$-mod}_\mathrm{tor}$ of $\mbox{\rm $\La(G)$-mod}.$ Thus one could define characteristic classes  for $\La(G)$-torsion modules inside this $K$-group. Unfortunately, in the
noncommutative case there are examples $M$ of (pseudo-null)
$\La(G)$-torsion modules whose class in
$K_0(\mbox{\rm $\La(G)$-mod}_\mathrm{tor})$ vanishes  though the
$G$-Euler-Poincar\'{e} characteristic of $M$ does not (cf.\
\cite[{\S}4]{cs2}), i.e.\ the latter $K$-group does not discover this characteristic and thus cannot bear the arithmetic content of  e.g.\
the Selmer group of an elliptic curve (see  section \ref{selmer}). 
In fact, the module in the following example which stems from \cite{cs2}  has this property. It also illustrates   that, in general, the $G$-Euler characteristic is not invariant under  pseudo-isomorphisms. 

\begin{ex}\label{ex-css}
Let $G=H\times \Gamma$  be a $p$-valuable group with $\Gamma\cong\zp,$ generated by $\gamma$ say, and assume that $H$ contains  a subgroup which is a semi-direct product of the following form. Let $H_1$ and $H_2$ be two closed subgroups of $H$ which are isomorphic to $\zp,$ and  which are such that $h_2h_1h_2^{-1}=h_1^{\phi(h_2)}$ for  fixed topological generators $h_i$ of $H_i,$ where $\xymatrix{
 {\phi:H_2\ar@{^{(}->}[r]} & {\mathrm{Aut}(H_1)=\mathbb{Z}_p^\times }}$ is a continuous injective group homomorphism; in particular,  the subgroup $H_1H_2$ of $H$ is non-abelian. For example, such subgroup exists for any $H$ which is open in $SL_n(\zp)$ $(n\geq 2.)$ Putting \[g=h_1-1,\; \omega:=h_2+p^r,\]
for any integer $r\geq 1,$ and \[u=h_2\cdot \frac{h_1^{\phi(h_2)}-1}{h_1-1}\] one verifies that \begin{eqnarray}\label{relation}
g\omega=ug
\end{eqnarray} holds. Using this relation one sees immediately that right multiplication by $g$ on $\La$ induces a short exact sequence 
\begin{eqnarray}\label{def-seq}
\xymatrix@1{
   {\ 0 \ar[r] } &  {\ \La/\La(\gamma-u) \ar[r]^{\cdot g} } &  {\ \La/\La(\gamma-\omega) \ar[r] } &  {\ M \ar[r] } &  {\ 0 } }
\end{eqnarray}
defining the $\La$-module $M.$ Note that $\La/\La(\gamma-\omega)$ is isomorphic as $\La(G)$-module to $\La(H)$ on which $\Gamma$ acts via right multiplication by $\omega,$ similar for $\La/\La(\gamma-u).$ Thus  $M$ is a pseudo-null $\La(G)$-module because its $\La(H)$-rank is zero, cf.\ \cite{ven-weier}. We postpone the calculation of the Euler characteristic to section \ref{desccent} where we will have  more techniques available, see Example \ref{ex-char} where we also construct a pseudo-null module with non-trivial Euler characteristic for the semi-direct product $\zp\rtimes\zp.$  
We determine the class of $M$ in $K_0(\mbox{\rm $\La(G)$-mod}_\mathrm{tor}):$ from the short exact sequence \eqref{def-seq} one sees  that the class of \[(\gamma-\omega)(\gamma-u)^{-1}\]
in $Q(G)^\times/[Q(G)^\times,Q(G)^\times]\La(G)^\times$ represents that module, and is zero. For the convenience of the reader we recall the short argument from \cite[{\S}4]{cs2}: Since $\gamma$ is in the center of $G,$ by \eqref{relation} we have \[g(\gamma-\omega)=(\gamma-u)g\] and thus we get in $Q(G)^\times/[Q(G)^\times,Q(G)^\times]\La(G)^\times$
\begin{eqnarray}
g(\gamma-\omega)(\gamma-u)^{-1}&=&(\gamma-u)g(\gamma-u)^{-1}\\
&=&(\gamma-u)(\gamma-u)^{-1}g\\
&=&g.
\end{eqnarray}
Since $g$ is invertible in $Q(G)^\times/[Q(G)^\times,Q(G)^\times]\La(G)^\times$ the result follows. But we should also mention that we do  not know any example for a pseudo-null module with non-vanishing class in $K_0(\mbox{\rm $\La(G)$-mod}).$
\end{ex}

Thus one could try to  search for smaller subcategories of $\mbox{\rm $\La(G)$-mod}_\mathrm{tor}.$ But  the characteristic class in the corresponding $K$-group can only be identified with an characteristic element if
the suitable subcategory of
$\mbox{$\La$-mod}_{\mathrm{tor}}$ can be described by
a pair of rings  as for $\La\to \La_\mathcal{T}$ for a suitable Ore-set $\T.$ Also, only in this case the formalism  of virtual objects can be applied to
``descent" as we will do in section \ref{desccent}. This is our motivation to study Ore-sets associated with certain group extensions in section \ref{Ore-sets}.
\end{rem}

\POP
\PUSH{relative.tex}%

\Section{Ore sets associated with group extensions}\label{Ore-sets}

In section \ref{swan} we saw that for every Ore-set $\T$ the group $K_0(\La,\La_\T)$ describes the Grothendieck group of $\T$-torsion \La-modules. From classical Iwasawa theory over $\zp$-extensions we know that characteristic elements live in $K_0(\La(\Gamma),Q(\Gamma)).$ In order to make use of this information we are looking for a ring $R$ with $\La(G)\subseteq R\subseteq Q(G)$ such that the projection $\psi_H:\La(G)\to\La(\Gamma)$ extends to a commutative diagram \[\xymatrix{
  {\La(G)\ar@{->>}[r]\ar@{^{(}->}[d]} & {\La(\Gamma)\ar@{^{(}->}[d]} \\
  R\ar[r] & Q(\Gamma). 
}\]

The first candidate for $R$ would be $\La_{\T'}$ in case $\T':=\La(G)\setminus \ker(\psi_H)$ satisfies the Ore condition. But firstly it seems difficult to prove this for a general class of groups (only the case $G=\zp\rtimes\zp$ is known and straightforward) and secondly even if the existence is known the associated $\T'$-torsion category is not closed under the kind of twisting by representations (even for $G$ abelian) that will be discussed  in section \ref{twisting}. Thus we shall take a slightly smaller set $\T$ below. 

Since these  type of questions  are of general interest and since there does not seem to exist any localisation result of this sort in the literature we treat this topic in greater generality than needed for our applications. But see \cite[thm.\ 2.14/15]{pass} where a similar topic is discussed in the context of the (usual) group algebra of polycyclic-by-finite groups with coefficients in a field.

Let $G$ be an extension of a torsionfree pro-$p$ $p$-adic Lie
group $\Gamma$ by a $p$-adic Lie group
 $H,$ i.e.\ we have an short exact sequence
\[
\xymatrix@1{
   {\ 1 \ar[r] } &  {\ H \ar[r] } &  {\ G \ar[r] } &  {\ \Gamma \ar[r] } &  {\
   1. }
}\]

The projection $G\twoheadrightarrow\Gamma$ induces  a canonical
surjective ring homomorphism \[\psi_H:\O\kl G\kr \twoheadrightarrow
\kappa\kl \Gamma\kr,\] where $\kappa$ denotes the residue class field of $\O$ as before. Note that due to compactness
$\mh:=\ker(\psi_H)$ equals $\O\kl G\kr\m_H=\m_H\O\kl G\kr,$ where
$\m_H$ denotes the kernel of the canonical map \[\O\kl H\kr
\twoheadrightarrow \kappa.\]

We put $\mathcal{T}:=\La\setminus\mh$ and   dare to formulate the following

\begin{conj}\label{Ore-conj} Assume that $H$ is a pro-$p$  group. Then, the multiplicative closed set $\mathcal{T}$ is an  Ore set of $\La.$
\end{conj}

Recall that an element $x\in R$ is {\em right
regular} if $xr=0$ implies $r=0$ for $r\in R.$ Similarly {\em left
regular} is defined and {\em regular} means both right and left
regular (and hence not a zero divisor). For an ideal $I$ of $R$ we
define the multiplicatively closed set
\begin{eqnarray*}
\mathcal{C}_R(I)&:=&\{s\in R| s+I \mbox{ is regular in } R/I\}.
\end{eqnarray*}
Note that the above set $\mathcal{T}$ is nothing else than
$\mathcal{C}(\mh)$ because $\kappa\kl \Gamma\kr\cong\kappa[[X]]$
is an integral domain or in other words  $\mh$ is a completely
prime ideal.

\begin{thm} \label{ore-thm} Let  $G$  the semi-direct product $\zp\rtimes\zp$ or the
direct product $H\times \Gamma$ of a uniform group $H$ with  a
torsion-free pro-$p$ $p$-adic Lie group $\Gamma.$ Then the
multiplicative set $\mathcal{T}$ satisfies the left and right Ore condition
for $\La=\O \kl G\kr$ for $\O$ either a finite field or a DVR. In
particular, the localisation $\La_\mathcal{T}$ exists and is a Noetherian
regular local ring with $\gl {\La_\mathcal{T}}<\gl \La.$
\end{thm}

\begin{proof}
We mention that since the uniformising element of $\O$ is central
in $\La$ and contained in $\mh$ it would suffice to prove the
statement for a finite field $\kappa$ by \cite[lem.\ 4.2]{smith}.
Now the strategy is to show that \begin{enumerate} \item $\mh$
satisfies the (left and right)  Artin Rees property
(\cite[4.2.2]{mc-rob}) and
\item $\mathcal{C}_{\La/\mh^2}(\mh/\mh^2)$ is an Ore set of
$\La/\mh^2.$
\end{enumerate}
Then \cite[cor.\ 4.7]{smith} (see also
\cite[4.2.10]{mc-rob}) implies that $\mathcal{T}$ is an Ore set of $\La.$ These properties will be investigated in the next subsections. The other statements follow from the following lemma.
\end{proof}

We hope that one can modify the above criteria (i) and  (ii) replacing the $\mh$-adic filtration by a filtration induced by Lazard's more general filtration associated with $\La(H)$ for any $p$-valued group $(H,\omega)$ or using the stronger criterion of Smith \cite[lem.\ 4.1, thm.\ 4.6]{smith} applied to the $\m(H')$-adic filtration associated with an uniform normal subgroup $H'$ of $G$ contained in $H.$ This hope is the reason for the above conjecture which could even extend to a larger class of not necessarily pro-$p$ groups, e.g.\ to open subgroups of $GL_n(\zp).$ In fact, there is already joint work with J. Coates and R. Sujatha in progress in order to settle these cases.

\begin{lem}\label{local}Let $\La$ be a Noetherian local ring with maximal ideal $\m$ and $\P\neq \m$ a completely prime ideal. Suppose that $\mathcal{T}:=\mathcal{C}(\P)=\La\setminus \P$ satisfies the Ore condition. Then  $\La_\mathcal{T}$ is again a  local ring, i.e.\ its non-units form a
maximal ideal, and the global dimension of $\La_\mathcal{T}$ is
(strictly) less than that of $\La.$ More precisely, the units are of the following form \[\La_\mathcal{T}^\times=\{\lambda t^{-1}|\; \lambda, t \in \mathcal{T}\}.\]
\end{lem}

\begin{proof} 
Let $\mathcal{M}\subsetneqq \La_\mathcal{T}$ be a maximal left ideal. By
\cite[Prop.\ 2.1.16]{mc-rob}  we conclude that
$\mathcal{M}=\La_\mathcal{T}(\mathcal{M}\cap \La)\subseteq I(H)_\mathcal{T}$ as $\mathcal{M}\cap
\La\subseteq \P$ by the properness of $\mathcal{M}$ (Note that this
argument holds only for localisations at completely prime ideals).
By symmetry we  see that $\P_\mathcal{T}$ is the unique left and
the unique right maximal ideal of $\La_\mathcal{T}$ which implies
that $\La_\mathcal{T}$ is local. The statement concerning the
global dimension follows from \cite[cor. 4.3]{mc-rob} and theorem 4.4
(loc.\ cit.) applied to the unique simple $\La$-module $\La/\m.$  
The description of the units is an immediate consequence from the well known fact that the maximal ideal of $\La_\mathcal{T}$ is $\P\La_\mathcal{T}=\{pt^{-1}|p\in \P, t\in \mathcal{T}\}.$
\end{proof}

For the next statement  we assume that $\Gamma\cong\zp.$  In the arithmetic applications we have in mind those $\La(G)$-modules which are finitely generated over $\La(H)$ play an important role, see Theorem \ref{tate-case}. We write $\mbox{\rm $\La(G)$-mod}^H$ for the full subcategory of \mbox{\rm $\La(G)$-mod} consisting of such modules. Recall that a $\La(G)$-module $M$  
is called $\mathcal{T}$-torsion, if   $M\otimes_{\La}\La_\mathcal{T}=0.$ The significance of $\mathcal{T}$ being an Ore set results from the following observation.

\begin{prop} \label{t-torsion} Assume that $\T$ is an Ore-set.
\begin{enumerate}
\item Let  $g\in M_n(\Lambda(G))\cap GL_n(\La(G)_{\T})$ for some  natural number $n,$
 Then the cokernel of $g$ is finitely generated as $\La(H)$-module while $\ker g$ is trivial.
 \item The inclusion $\La(G)\mbox{-mod}^H\subseteq
 \La(G)\mbox{-mod}_{\T-\mathrm{tor}}$ is
 an identity of categories where the latter one consists of all
 those
 finitely generated $\La(G)$-modules which are
 $\La(G)_{\T}$-torsion.
\item $K_0(\Lambda(G),\La(G)_{\T})\cong
K_0(\La(G)\mbox{-mod}_{\T-\mathrm{tor}})\cong
K_0(\La(G)\mbox{-mod}^H).$
\end{enumerate}
\end{prop}

Before we give the proof consider the canonical ring homomorphism
$\pi_H:\La(G)\to\La(\Gamma)$ which is induced by the group
homomorphism $G\to\Gamma.$ Since $\pi_H^{-1}((p))=\mh$ and thus $\pi_H(\T)\subseteq \La(\Gamma)\setminus (p),$ one obtains a ring homomorphism \[\La(G)_{\T}\to\La(\Gamma)_{(p)}.\]

\begin{proof}
Set $M:=\coker(g).$ For all $n$  the augmentation map $\pi_H$
induces homomorphisms of groups
\[GL_n(\La(G)_{\T})\to GL_n(\La(\Gamma)_{(p)}).\]
Thus reduction modulo $I(H)$ induces a short exact sequence
\[\xymatrix@1{
   {\ 0 \ar[r] } &  {\ \La(\Gamma)^n \ar[r] } &  {\ \La(\Gamma)^n  \ar[r] } &  {\ M_H \ar[r] } &  {\ 0
   },
}\] where $M_H$ is a finitely generated $\zp$-module because any
element in $M_H$ has an annihilator prime to $p$ (this implies the left-exactness because  $M_H$ being  $\La(\Gamma)$-torsion forces the kernel to be torsion as well, but $\La(\Gamma)^n $ is torsionfree). Then Nakayama's
lemma implies (i) (for the kernel use the same argument as above).   Next we prove that every module in $\La(G)\mbox{-mod}^H$ is indeed
$\La(G)_{\T}$-torsion:  Let $m\in M$ arbitrary. By Lemma  \ref{cyclic-wp} below
 there is an element in $\T$   and a surjection $\La/\La f \twoheadrightarrow \La
m \subseteq M.$ Since $f$ obviously annihilates $m$ the claim follows.

Now let $M\in \La(G)\mbox{-mod}_{\T-\mathrm{tor}}$ be arbitrary and choose a finite set of $\La(G)$-generators $m_i, i\in I$ of $M.$ By assumption there exists for every $i\in I$ an element $f_i$ in $\T$ such that $f_i m_i=0$ and thus $M$ is the homomorphic image of the finitely generated $\La(H)$-module $\bigoplus_{i\in I} \La(G)/\La(G)f_i$ by Lemma \ref{cyclic-wp}   again. This proves (ii). The last item follows from (ii) and the standard
exact localisation sequences of (relative) $K$-theory.
\end{proof}

\begin{lem}\label{cyclic-wp}
Let $J$ be a left ideal of $\La=\La(G).$ Then $M:=\La/J$ is
finitely generated as $\La(H)$-module if and only if $J$ contains an element of $\T.$
 \end{lem}
 
 \begin{proof}
 Asssume first that $M$ is a finitely generated $\La(H)$-module and suppose 
   all elements of $J$ reduce to zero in $\kappa\kl \Gamma\kr.$ Then $M/\m(H) M\cong\kappa\kl \Gamma\kr.$   But since there is some surjection $\La(H)^n\to M$ and since $\m_H\subseteq\m(H)$ the module $M/\m(H) M$ is a finitely generated $\kappa$-module, a contradiction. 
Now assume for simplicity that $J=\La(G)f$ with $f\in\T,$ i.e.\ we have a short exact sequence \[\xymatrix@1{
    {\ \La(G) \ar[r]^{\cdot f} } &  {\ \La(G)\ar[r] } &  {\ M\ar[r] } &  {\ 0}. 
}\]
Tensoring with $\kappa\kl \Gamma\kr\cong \La(G)/\m(H)$ leads to 
\[\xymatrix@1{
    {\ \kappa\kl \Gamma\kr \ar[r]^{\cdot \bar{f}}  } &  {\ \kappa\kl \Gamma\kr \ar[r]} &  {\ M/\m_H M \ar[r] } &  {\ 0}
}\]
showing that $M/\m_H M$ is a finite-dimensional $\kappa$-module. Thus, by Nakayama's Lemma, $M$ is finitely generated over $\La(H).$
 \end{proof}

\Subsection{The direct product case}

\begin{thm}\label{productcase} Let  $G=H\times\Gamma$ be the direct product of a
uniform group $H$ and a torsion-free pro-$p$ $p$-adic Lie group
$\Gamma.$ Then the Iwasawa algebra  $\La(G)$ is complete and
separated with respect to its $\mh$-adic filtration. Its
associated graded ring $\gr{\La(G)}$ is isomorphic to a
generalized enveloping algebra
\[\gr_{\mh}{\La(G)}\cong \left\{ \begin{array}{ll}
U(\gr H \otimes_{\gr\zp}\kappa[\pi])\otimes_{\kappa[\pi]}\kappa\kl \Gamma \kr[\pi]& \mbox{ if $\O$ is a DVR,} \\
U(\gr H \otimes_{\gr\zp}\kappa)\otimes_{\kappa}\kappa\kl \Gamma
\kr & \mbox{
 if $\O=\kappa$ is a finite field.}
\end{array}\right.
\]

 in
particular, it is a Noetherian integral domain.
\end{thm}


For $G=H\rtimes\Gamma,$ the above isomorphism is still an isomorphism of grade $\kappa[\pi]$-modules. Analyzing the induced ring structure will hopefully lead to a proof of Conjecture \ref{Ore-conj} in this case.

\begin{proof}
Since $\mh$ is contained in the maximal ideal of $\La(G)$ the
$\mh$-adic filtration is separated. By compactness of $\La(G)$ the
canonical map \[\La(G)\to\projlim n \La(G)/\mh^n\] is thus an
isomorphism. Now we calculate the graded ring, using flatness of
$\La(G)$ over $\La(H)$   and Corollary \ref{uniform} we obtain the following
isomorphisms of graded $\kappa[\pi]$-modules
\begin{eqnarray*}
\gr_{\mh}{\La(G)}&\cong&\bigoplus_{n\geq 0} \mh^n/\mh^{n+1}\\
&\cong&\bigoplus_{n\geq 0} \m_H^n\La(G)/\m_H^{n+1} \La(G)\\
&\cong&\bigoplus_{n\geq 0} \Big((\m_H^n/\m_H^{n+1})
\otimes_{\La(H)}\La(G)\Big)\\
&\cong&\bigoplus_{n\geq 0} \Big((\m_H^n/\m_H^{n+1})
\otimes_{\kappa}\kappa\kl \Gamma\kr \Big)\\
&\cong& \Big(\gr_{\m_H}\La(H)\Big)\otimes_\kappa \kappa\kl
\Gamma\kr \\
&\cong& U(\gr H\otimes_{\gr\zp}\gr \O)\otimes_\kappa
\kappa\kl\Gamma\kr \\
&\cong& U(\gr H\otimes_{\gr\zp}\gr \O)\otimes_{\kappa[\pi]}
\kappa\kl\Gamma\kr[\pi].
 \end{eqnarray*}

Since $G=H\times \Gamma$ this is also a ring-isomorphism. By
\cite[1.7.14]{mc-rob} and \ref{lazard} the latter ring is an
Noetherian integral domain. Observe that the above calculation   holds (up to the ring structure) also for arbitrary group extensions of the form  considered above.
\end{proof}

Recall that a filtered ring $R$  with filtration $F_\bullet R$ is
called a Zariski ring if its associated Rees ring \[
\widetilde{R}=\bigoplus_{n\in\mathbb{Z}} F_nR t^n\subseteq
R[t,t^{-1}]\] is Noetherian and $F_{-1}R$ is contained in the
Jacobson ideal of $F_0R.$ Now, the equivalences of (1), (3) and
(4) of \cite[ch.\ II {\S}2 thm.\ 2.1.2]{li-o} imply the following

\begin{cor}\label{zariski}
 Under the assumptions of the theorem $\La(G)$ endowed with its $\mh$-adic filtration is a
Zariski ring, in particular it satisfies the (left and right)
Artin Rees property for the ideal $\mh.$
\end{cor}

\begin{lem}\label{lem-Ore}
Under the conditions of the theorem it holds that
\[\mathcal{C}_\La(\mh)\subseteq \mathcal{C}_\La(\mh^2).\]
Thus $\mathcal{C}_{\La/\mh^2}(\mh/\mh^2)$ is an Ore set of
$\La/\mh^2.$
\end{lem}

\begin{proof}
Let $\lambda$ be in $\mathcal{C}_\La(\mh)=\La\setminus\mh$ and let
$\lambda'$ be an element in $\La$ such for which
$\lambda'\lambda\in\mh^2$ holds.  Note that $\lambda'$ is in $\mh$
because $\lambda$ belongs to $\mathcal{C}_\La(\mh).$ We have to
prove that $\lambda'$ is  in $\mh^2.$ We assume the contrary,
i.e.\ $\lambda'\in\mh\setminus\mh^2.$ But then we obtain that in the
graded ring $\gr_{\mh}\La(G)$
\[(\lambda'+\mh^2)\cdot (\lambda +\mh)=\lambda'\lambda+\mh^2=0,\]
which contradicts the integrality of that ring.

The implication that if $\lambda\lambda'\in\mh$ then
$\lambda'\in\mh$ follows by symmetry and thus we have shown the
first statement which in turn implies that
\[\mathcal{C}_{\La/\mh^2}(\mh/\mh^2)\subseteq \mathcal{C}_{\La/\mh^2}(0).\]
Now the second statement follows from Small's theorem
\cite[4.1.3/4]{mc-rob}.
\end{proof}

\Subsection{The semi-direct product case}

In this subsection we restrict to the easiest semi-direct product
case, viz $G=\zp\rtimes\zp,$ though the methods certainly extend
to a wider class of poly-cyclic pro-$p$-groups.

Again, we  are first concerned with the Artin-Rees property.

\begin{prop}
Let $\O=\kappa$ be a finite field. Then $\mh$ satisfies the Artin
Rees property.
\end{prop}

\begin{proof}
We identify $\La(G)$ with the skew power series ring
$\kappa[[Y,X;\sigma,\delta]]$as in \cite{ven-weier}. Note that $Y$ is a normal element
of $\La(G),$ i.e.\ $\La(G)Y=Y\La(G),$ which generates $\mh.$ Thus the statement follows from
\cite[thm.\ 4.2.7]{mc-rob}.
\end{proof}

With the same technique and some calculations one easily shows
that the proposition holds also if $\O$ is a DVR.

\begin{lem}
For $\O$ a DVR or a finite field and a uniform group
$G=H\rtimes\Gamma$ which is the semidirect product of  a normal
uniform subgroup $H$ and $\Gamma\cong\zp$  it holds that
\[\mathcal{C}_\La(\mh)\subseteq \mathcal{C}_\La(\mh^2).\]
Thus $\mathcal{C}_{\La/\mh^2}(\mh/\mh^2)$ is an Ore set of
$\La/\mh^2.$
\end{lem}

\begin{proof}
Let $\lambda$ be in $\mathcal{C}_\La(\mh)=\La\setminus\mh$ and let
$\lambda'$ be an element in $\La$ such that
$\lambda'\lambda\in\mh^2$ holds.  Note that $\lambda'$ is in $\mh$
because $\lambda$ belongs to $\mathcal{C}_\La(\mh).$ We have to
prove that $\lambda'$ is  in $\mh^2.$ We assume the contrary,
i.e.\ $\lambda'\in\mh\setminus\mh^2.$
We identify $\La(G)$ with the skew power series ring
$\La(H)[[X;\sigma,\delta]]$ (see \cite{ven-weier}, $\sigma$ is the ring isomorphism of $\La(H)$ induced by the operation of $\Gamma$ on $H$ while $\delta$ denotes the $\sigma$-derivation $\sigma -\mathrm{id}$). We expand $\lambda$ and $\lambda'$ as
\[\lambda=\sum_{i\geq 0} \lambda_i X^i \mbox{ and }\lambda'=\sum_{i\geq 0} \lambda'_i
X^i,\] where all $\lambda'_i\in\m_H$ (note that for any $n\geq 0$
the ideal $\mh^n$ consists precisely of those power series in the
variable $X$ whose coefficients all lie in $\m_H^n$). By
assumption there exist $i_0$ and $j_0$ such that $\lambda_{i_0}$
and $\lambda'_{j_0}$ are not in $\m_H$ and $\m_H^2,$ respectively.
Let us assume that these indices are chosen minimal with this
property. We want to calculate the product $\lambda'\lambda$ in
$\La/\mh^2$ and we observe that the latter ring is isomorphic to
the skew power series ring $\La(H)/\m_H^2[[X;\overline{\sigma}]]$
where $\overline{\sigma}$ is induced by $\sigma,$ while $\delta$
induces the zero derivation. For the $(i_0+j_0)$th coefficient one
obtains
\begin{eqnarray*}
0\equiv(\lambda'\lambda)_{i_0+j_0}&\equiv&\sum_{k+l=i_0+j_0}
\lambda'_k\sigma^k(\lambda_l) \mbox{ mod } \m_H^2.
\end{eqnarray*}
The products
$\lambda'_k\sigma^k(\lambda_l)$ are in $\m_H^2$ for $k<j_0$ by definition of $j_0$ and
 for $k>j_0$ because then $\lambda'_k$ and
$\sigma^k(\lambda_l)$ ($l<i_0$) both belong to $\m_H$ as $\m_H$ is
$\sigma$-invariant. Thus
$\lambda'_{j_0}\sigma^{j_0}(\lambda_{i_0})$ belongs to $\m_H^2,$
which is a contradiction as $\sigma^{j_0}(\lambda_{i_0})$ is a
unit of $\La(H).$ The rest is identical as in the proof of lemma
\ref{lem-Ore}.
\end{proof}

%
\POP

\PUSH{twisting.tex}%
\Section{Twisting} \label{twisting}

The twisting of the complex $L$-function by an Artin character corresponds on the algebraic side to tensoring the associated modules (e.g.\ the dual of the Selmer group of an elliptic curve) by the accordant representation. Basic properties of the latter formalism are studied in the first  subsection. In the second subsection we apply it in order to define the evaluation of an characteristic class or element at certain $p$-adic representations. In the third subsection we discuss several definitions of (equivariant) Euler-characteristics.

\Subsection{Twisting of \La-modules}

Let $T$ be a free $\O$-module of finite rank $r$ with continuous $G$-action given by 
\[\rho:G\rightarrow \mathrm{Aut}_\O(T)=Gl_r(\O).\]

 \begin{defn}
 For a finitely generated $\La=\La(  G)$-module $M$ we define the finitely generated
 \La-module $$M(\rho):=M\otimes_{\O}T=\Hom_{cont.,\O}(M,A)^{\vee}$$ with diagonal
 $  G$-action and where $A=T^\vee$ denotes the Pontryagin dual of $T,$
 \end{defn}

 Note that the functor $-(\rho)$ is exact. The following lemma is well-known.

\begin{lem} \label{twist-free} For any choice of an $\O$-basis of $T,$ i.e.\ of an $\O$-module isomorphism $\phi: T\cong\O^r,$ there is a canonical isomorphism \[\Lambda(\rho)\to \La\otimes_\O \O^r\cong\La^r,\] induced by mapping $g\otimes t,$ $g\in G,\; t\in T,$ to $g\phi(\rho(g^{-1})t).$
  \end{lem}

%
%

From this lemma, it follows that if $P$ is a projective
$\Lambda$-module, then so is $P(\rho)$.

Now let $\rho: G\to GL(V)$ be a continuous linear representation on a
finite dimensional Vector space $V$ over a finite extension $K$ of
$\qp$ with ring of integers $\O.$ We choose a $G$-invariant
$\O$-lattice $T\subseteq V.$ We denote by $\La(G)=\O\kl G\kr$ the
Iwasawa-algebra of $G$ with coefficients in $\O.$

The following lemma is immediately verified.

\begin{lem}For any finitely generated $\La$-module $M$ there are
canonical isomorphism of $\O$-modules
  \begin{eqnarray} T\otimes_{\O\kl G \kr} M&\cong&(T\otimes_\O M)\otimes_{\mathcal{O}\kl G\kr} \O,
\end{eqnarray}
where $t\otimes m$ is mapped to $t\otimes m\otimes 1$ while the inverse map  is induced by mapping $t\otimes m\otimes o$ to $o(t\otimes m)=(ot)\otimes m= t\otimes (om).$ They induce
isomorphisms
 \begin{eqnarray}
 \Tor_i^{\O\kl G \kr}(T,M)&\cong&\Tor_i^{\mathcal{O}\kl G \kr}(T\otimes_\O M,\O)
\end{eqnarray}
for all $i\geq 0.$
\end{lem}

Now let us again assume to be in the    situation of section \ref{Ore-sets} with $G$ being pro-$p$ and $G/H\cong\zp.$
The exact functor $T\otimes_\O -$ induces the following homomorphism of
  $K$-groups
  \[\rho_\ast:K_0(\La(G)\mbox{-mod}^H)\to K_0(\La(G)\mbox{-mod}^H),\]
indeed, if $M$ is finitely generated over $\La(H),$ so is $T\otimes_\O M.$ 
In fact, $\rho_\ast$ is independent of the choice of the lattice $T.$ 

\begin{lem}\label{lattice-indep}
Assume that the Weierstrass preparation theorem holds or that $\kappa\kl G\kr$ is an integral domain. Let $M\in\La(G)\mbox{\rm -mod}^H$ and $T'$ a further $G$-invariant lattice of $V.$ Then we have \[[T\otimes_\O M]=[T'\otimes_\O M]\] in $K_0(\La(G)\mbox{\rm -mod}^H).$
\end{lem}

\begin{proof}
By a standard argument we may assume that $T'\subseteq T$ and that $M=\La/\La f$ with an element $f$ of finite reduced order because classes of such modules generate the $K$-group. 
It follows easily from our assumption that $M$ is a torsionfree $\O$-module. 
Denoting by $E$ the finite $G$-module $T/T'$ we obtain thus the following exact sequence 
\[\xymatrix@1{
   {\ 0 \ar[r] } &  {\ T'\otimes_\O M \ar[r] } &  {\ T\otimes_\O M  \ar[r] } &  {\ E\otimes_\O M \ar[r] } &  {\ 0. } 
}\]
Using a Jordan H{\"o}lder series of $E$ it is sufficient to show that $[\kappa\otimes_\O M]$ vanishes. But this follows from the following exact sequence in $\La(G)\mbox{-mod}^H$
\[\xymatrix@1{
   {\ 0 \ar[r] } &  {\ M \ar[r]^{\pi} } &  {\ M \ar[r] } &  {\ \kappa\otimes_\O M \ar[r] } &  {\ 0, } 
}\]
where $\pi$ denotes an uniformizer of $\O.$
\end{proof}

Using the isomorphism $K_0(\Lambda(G),\La(G)_{\T})\cong
K_0(\La(G)\mbox{-mod}^H)$ established in Pro\-position \ref{t-torsion}, we obtain also a  homomorphism \[\rho_\ast:K_0(\Lambda(G),\La(G)_{\T})\to K_0(\Lambda(G),\La(G)_{\T}),\] which can  be described as follows:

Consider a triple  $(P_1,P_2,\lambda)$ representing a class of $K_0(\Lambda(G),\La(G)_{\T})$ with $P_i$ projective (thus free) $\La(G)$-modules and $\lambda:P_1\otimes_{\La(G)}\La(G)_{\T}\to P_2\otimes_{\La(G)}\La(G)_{\T}$ an isomorphism of
$\La(G)_{\T}$-modules.  Since $P_i\otimes_{\La(G)}\La(G)_{\T},\; i=1,2,$ are free $\La(G)_{\T}$-modules of rank $m,$ say, and since  $\lambda$ can be described by a invertible matrix with coefficients in $\La(G)_{\T}$ it is easily seen by finding a common denominator of the matrix elements that there exist matrices $A_i\in M_m(\La(G))\cap GL_m(\La(G)_{\T})$  such that $A_1(A_2)^{-1}$ represents $\lambda$ (for a certain choice of bases). Now we twist the $\La(G)$-homomorphisms $\lambda_i$ given by the $A_i$ with $T$ and denote the composite $T\otimes_{\O}\lambda_1 \circ (T\otimes_{\O}\lambda_2)^{-1}$ by   $T\otimes_\mathcal{O}\lambda.$ Now the triple $(P_1,P_2,\lambda)$ is sent to
$(T\otimes_\mathcal{O}P_1,T\otimes_\mathcal{O}P_2,T\otimes_\mathcal{O}\lambda).$

We also would like to see how  twisting  by $\rho$ operates on $K_1(\La(G)_{\mathcal{T}})\cong(\La(G)_{\mathcal{T}}^\times)^{ab} .$ To this end, for a finite dimensional continuous $\O$-representation $\rho:G\to
\mathrm{Aut}_O(T)$ we define the twist operator
\[\mathrm{tw}_\rho:\La(G) \to \mathrm{End}_\mathcal{O}(T)\otimes_\O \La(G)\]
as follows. By continuity we may assume that $\lambda\in\La(G)
  $ is of the form $\sum a_g g$ where almost all $a_g\in\O$ are
  zero. Then we set $\mathrm{tw}_\rho(\lambda):= \sum a_g\rho(g^{-1})\otimes g.$

The restriction of $\mathrm{tw}_\rho$ to
$\T$
  and the choice of an $\O$-basis of $T$ induces the multiplicative map
\[\mathrm{tw}_\rho:\T \to M_m(\La)\cap GL_m(\La_\mathcal{T}), \]
if $\mathcal{T}$ is an Ore-set of $\O\kl G\kr$ and where
$m=\rk_\mathcal{O}T.$ Thus $\mathrm{tw}_\rho$ extends to a ring homomorphism \[\mathrm{tw}_\rho:\La_\mathcal{T} \to M_m(\La_\mathcal{T}). \]
Restricting it to the units we obtain a group homomorphism
\[\mathrm{tw}_\rho:\La_\mathcal{T}^\times \to GL_m(\La_\mathcal{T}).\]
If we compose this map with the ``determinant" \[ \det:GL_m(\La_\mathcal{T})\to K_1(\La_\mathcal{T})\to (\La_\mathcal{T}^\times)^{ab}\] we have 
\[\det\circ\mathrm{tw}_\rho:\La_\mathcal{T}^\times \to(\La_\mathcal{T}^\times)^{ab},\]
which explicitly describes the action on $K_1(\La_\mathcal{T}).$ From a   functorial point of view this action is nothing else that the $K_1(-)$-functor applied to the ring homomorphism $tw_\rho$ and using Morita equivalence: \[K_1(\La_\mathcal{T}) \to K_1(M_m(\La_\mathcal{T}))\cong K_1(\La_\mathcal{T}).\]

In fact, we can easily extend componentwise
the above twisting operator to
\[\mathrm{tw}_\rho:M_n(\La(G)) \to \mathrm{End}_\mathcal{O}(T)\otimes_\O
M_n(\La(G)).\]
Note that the latter ring can be identified - for a chosen
$\O$-basis of  $T$ - with $M_{nm}(\O\kl G\kr).$ The augmentation
map $\pi_H: \O\kl G\kr\to \O\kl \Gamma\kr$ induces the map
\[\pi_H: M_n(\O\kl G\kr)\to M_n(\O\kl \Gamma\kr)\] which we denote
by the same symbol by abuse of notation.

Then we have the more general

\begin{lem}\label{twisting-lem}
 Let $g$ be in $M_n(\La(G))$ for some $n$ and assume
that $\coker(g)$ is finitely generated as $\La(H)$-module. Then
$\mathrm{tw}_\rho(g)$ is in $GL_{nm}(\La_\mathcal{T})$  if $\mathcal{T}$   is an Ore-set of $\La.$ In any
case \[\pi_H(\mathrm{tw}_\rho(g))\in
GL_{nm}(Q(\Gamma)).\]
\end{lem}

\begin{proof}
Applying $T\otimes_\O -$ to the short exact sequence
\[\xymatrix@1{
   {\ 0 \ar[r] } &  {\ \La(G)^n \ar[r]^{g} } &  {\ \La(G)^n \ar[r] } &  {\ \coker(g) \ar[r] } &  {\ 0. }
   }\]
gives the short exact sequence
\[\xymatrix@1{
   {\ 0 \ar[r] } &  {\ \O\kl  G \kr^{nm} \ar[r]^{\mathrm{tw}_\rho(g)} } &  {\ \O\kl  G \kr^{nm} \ar[r] } &  {\ T\otimes_\O \coker(g) \ar[r] } &  {\ 0. }
   }\]
after choosing a $\O$-basis of $T$ and using the isomorphism
\[T\otimes_\O\La(G)\cong \O\kl \Gamma\kr^m\] which is induced by
$t\otimes g\mapsto \rho(g^{-1})t\otimes g$ with $g\in G$ and $t\in
T$ (cf.\ \ref{twist-free}). It is well-known that under our
assumptions $T\otimes_\O \coker(g)$ is again finitely generated
as $\La(H)$-module which implies the first statement. Taking
$H$-coinvariants we obtain the short exact sequence
\[\xymatrix@1@+20pt{
   {\ 0 \ar[r] } &  {\ \O\kl \Gamma\kr^{nm} \ar[r]^{\pi_H(\mathrm{tw}_\rho(g))} } &  {\ \O\kl \Gamma\kr^{nm} \ar[r] } &  {\ (T\otimes_\O \coker(g))_H \ar[r] } &  {\ 0. }
   }\]
which is injective since the kernel is $\O\kl
\Gamma\kr$-torsionfree of rank zero. Since $(T\otimes_\O
\coker(g))_H$ is a finitely generated $\O$-module the last claim
follows after tensoring with $Q_\mathcal{O}(\Gamma).$
\end{proof}

%

\begin{lem}\label{twchar}
Let $M$ be in $\La\mbox{-mod}^H$ and $T$ as above. Then, for every choice of characteristic elements $F_{T\otimes_\O M}$ and $F_M,$ we have the equality \[F_{T\otimes_\O M}=\mathrm{det}_{\La_\mathcal{T}}\circ\mathrm{tw}_\rho 
(F_M)\] in $(\La_\mathcal{T}^\times)^{ab}/\mathrm{im}(\La^\times).$
\end{lem}

\begin{proof}
By multiplicativity and additivity  of $F_M$ and $[M],$ respectively, we may assume that $F_M\in\mathcal{T}.$ Then $[M]=[\coker(F_M)]$ and applying $\rho_\ast$ gives immediately $[T\otimes_\O M ]=[\coker(\mathrm{tw}_\rho(F_M)]$ using that  for every $g\in\mathcal{T}$ \[T\otimes_\O \coker(g)\cong\coker(\mathrm{tw}_\rho(g))\] via the isomorphism \ref{twist-free}. Now the result follows from the definition of the determinant.
\end{proof}

\POP
\PUSH{evaluating.tex}%
\Subsection{Evaluating at representations} \label{evaluation}

We would like to evaluate the characteristic elements
$F_M\in\La(G)_\mathcal{T}$ at ``$0$" - in other words at the
trivial representation -, i.e.\ we would like to extend the
augmentation map $\La(G)\to\O$ to $\La(G)_\mathcal{T}\to K,$ in
order to relate the value to the $G$-Euler characteristic of $M.$
But since we have to be careful with those denominators which map
to zero we define it in two steps, firstly for elements of $Q(\Gamma)$ and secondly of $\La(G)_\mathcal{T}:$

First let $F=hg^{-1}$ be an element of $Q(\Gamma)$ with $h$ and $g$ in $\La(\Gamma)$ prime to each other. If the image $g(0)$ of $g$ under the augmentation $\La(\Gamma)\to \O$ is not zero, we say that $F$ can be evaluated at zero and set $F(0):=h(0)g(0)^{-1} .$ With other words, $F(0)$ is defined, if $F$ belongs to the localisation $\La(\Gamma)_{I(\Gamma)}$ at the augmentation ideal $I(\Gamma),$ and then equals the image of $F$ under the extended augmentation map $\La(\Gamma)_{I(\Gamma)}\to K.$

Now let $F=hg^{-1}$ be an element of $\La(G)_\mathcal{T}$ with $h$  in $\La(G)$ and $g$ in $\mathcal{T}\subseteq\La(G).$ We say that $F(0)$ is defined if $\pi_H(F)(0)$ is defined and then we set $F(0):=\pi_H(F)(0),$ where \[\pi_H:\La(G)_\mathcal{T}\to Q(\Gamma)\] is the extended augmentation map (with respect to $H$).

We want to evaluate $F\in\La_\mathcal{T}$ at certain representations using the twist operator. Thus let $\rho: G\to GL(V)$ be a continuous linear representation on a
finite dimensional vector space $V$ of dimension $m$ over a finite extension $K$ of
$\qp$ with ring of integers $\O.$ We choose a $G$-invariant
$\O$-lattice $T\subseteq V$ and fix an $\O$-basis.

Consider the following composition of ring homomorphisms
\[ \xymatrix@1{
   {\ \La_\mathcal{T} \ar[r]^{\mathrm{tw}_\rho }} &  {\ M_m(\La_\mathcal{T}) \ar[r]^{\pi_H}} & { M_m(Q(\Gamma))}. }\]
It induces a group homomorphism
\[ \xymatrix@1{
   {\ \La_\mathcal{T}^\times \ar[r]^-{\mathrm{tw}_\rho }} &  {\ GL_m(\La_\mathcal{T}) \ar[r]^-{\pi_H}} & { GL_m(Q(\Gamma))\ar[r]^-{\mathrm{det}}} & Q(\Gamma)^\times,}\]
which  factorizes over the abelianization $(\La_\mathcal{T}^\times)^{ab}$ of   $\La_\mathcal{T}^\times.$   Moreover, the composite $\det_{Q(\Gamma)}\circ\pi_H\circ\mathrm{tw}_\rho$ factorizes $\mbox{mod} \; \La(\Gamma)^\times$ over the quotient $(\La_\mathcal{T}^\times)^{ab}/\mathrm{im}(\La^\times).$ We say that $F(\rho)$ is defined if $(\det_{Q(\Gamma)}\circ\pi_H\circ\mathrm{tw}_\rho (F))(0)\in Q(\Gamma)$ is so, and then we take the latter as value of $F$ at $\rho.$ 
Since the change to another lattice and $\O$-basis corresponds to conjugation by some matrix $A\in GL_2(\qp)$ and due to taking determinant, the value of $F(\rho)$ in $K^\times$ is independent of the  choice of $T$ and a basis. 

From the functoriality of determinants and by Lemma \ref{twchar} we obtain immediately the following

\begin{lem}
\begin{enumerate}
\item For every $F$ in $\La_\mathcal{T}^\times$ it holds that \[\pi_H\circ\mathrm{det}_{\La_\mathcal{T}}\circ\mathrm{tw}_\rho(F)=\mathrm{det}_{Q(\Gamma)}\circ\pi_H\circ\mathrm{tw}_\rho (F).\]
\item For every $M\in \La\mbox{-mod}^H$ and $T$ as above we have, for   every choice of characteristic elements $F_{T\otimes_\O M}$ and $F_M,$ \[F_M(\rho)=F_{T\otimes_\O M}(0)\] in $K^\times/O^\times.$
\end{enumerate}
\end{lem}

There is a second but less general way to describe the evaluation at representations.
The continuous homomorphism $\rho:G\to GL_{\mathcal{O}}(T)$ of
groups induces a homomorphism \[\rho:\La(G)\to
\mathrm{End}_{\mathcal{O}}(T)\] of $\O$-algebras and we can
compose this map with
\[\mathrm{det}_K:\mathrm{End}_{\mathcal{O}}(T)\to \O\subseteq K \subseteq \overline{\qp}\]
in order to evaluate elements $f\in\La(G)$ at the representation
$\rho:$ \[f(\rho):=\mathrm{det}_K(\rho^d(f))\in
\O_{\overline{\qp}},\] where $\O_{\overline{\qp}}$ denotes the
ring of integers of $\overline{\qp}$ and where we have chosen the
contragredient representation $\rho^d$ of $\rho,$ which belongs to
the  representation on $\Hom_\mathcal{O}(T,\O)$ induced by $\rho,$
i.e.\ in terms of matrices $\rho^d(g)$ is the transpose matrix of
$\rho(g^{-1})$ (after choosing a basis of $T$). The reason for
this is that we want   compatibility  with the twisting
operator, see Lemma \ref{compa}. Note that this
value is independent of the choice of the lattice $T.$ Also one
easily verifies that \[f(\rho^d)=f^\iota(\rho)\] where
$-^\iota:\La(G)\to\La(G)$ denotes the involution  which maps $g$
to $g^{-1}.$

We still have to justify that this new definition  is
compatible with our earlier definition of $f(\rho).$

\begin{lem}\label{compa}
Let $F=fg^{-1}$ be in $\La_\mathcal{T}$ with $f\in\La(G),g\in\T.$ If $g(\rho)\neq 0,$
then we have \[F(\rho)=\frac{f(\rho)}{g(\rho)}.\] In particular, the quotient is independent of the choice of fraction $F=fg^{-1}$ with $g(\rho)\neq 0$ (if such fraction exists at all).
\end{lem}

\begin{proof}
If $g(\rho)\neq 0,$ then we have
\begin{eqnarray*}
\pi_\Gamma
(\mathrm{det}_{Q(\Gamma)}(\pi_H(\mathrm{tw}_\rho(F))))&=&
\pi_\Gamma
(\mathrm{det}_{Q(\Gamma)}(\pi_H(\mathrm{tw}_\rho(f))\cdot
\pi_H(\mathrm{tw}_\rho(g))^{-1}))\\&=&\pi_\Gamma
(\mathrm{det}_{Q(\Gamma)}(\pi_H(\mathrm{tw}_\rho(f))))\cdot
\pi_\Gamma (\mathrm{det}_{Q(\Gamma)}(
\pi_H(\mathrm{tw}_\rho(g))))^{-1}\\
&=&\mathrm{det}_K (\pi_G(\mathrm{tw}_\rho(f)))\cdot \mathrm{det}_K
(\pi_G(\mathrm{tw}_\rho(g)))^{-1}\\
&=&\frac{f(\rho)}{g(\rho)}.
\end{eqnarray*}
Note that the last equality relies on our choice to define
$f(\rho)$ using the contragredient representation $\rho^d$ and
thus we have $\pi_G(\mathrm{tw}_\rho(f))=\rho^d(f)^t\in
\mathrm{End}_\mathcal{O}(T)$ where $t$ indicates the dual or
transpose endomorphism.
\end{proof}

\POP
\PUSH{ec.tex}%
\Subsection{Equivariant Euler-characteristics} \label{eq-EC}

Let $\O$ be a complete discrete valuation ring which is a finitely
generated   \zp-module, and let $K$ be its field of quotients. Let $G$
a compact (not necessarily pro-$p$) $p$-adic Lie group and $U\subseteq G$ a normal open
subgroup. Then we set $\Delta:=G/U.$ It is well known as a
consequence of Maschke's theorem and general Wedderburn theory
that the group algebra with coefficients in $K$ decomposes as a
product of matrix algebras
\[K[\Delta]\cong \prod_{l=1}^k M_{n_l}(D_l^{o}),\]
where $D_l^{o}$ denotes the   opposite of the division fields of
endomorphisms $\mathrm{End}_{K[\Delta]}(V_l)$  corresponding to a
system of representatives of irreducible $K$-representations
$V_l,$  $1\leq l\leq k,$ of $\Delta.$  The integers $n_l$ is the
length of $\mathrm{End}_{D_l}(V_l)$ and equals the multiplicity
with which the representation $V_l$ occurs in the regular
representation of $K[\Delta].$ The centers $K_l$ of
$D_l=\mathrm{End}_{K[\Delta]}(V_l)$ (and
$M_{n_l}(D_l^{o})=\mathrm{End}_{D_l}(V_l)$) are finite field
extensions of $K,$ whose ring of integers we denote by $\O_l.$

Henceforth we suppose that $D_l=K_l$ for all $1\leq l\leq k$ which
holds e.g.\ if $K$ is a splitting field of $\Delta$ but we don't
assume this stronger condition. Then $\prod_{l=1}^k M_{n_l}(\O_l)$ is a maximal order in $K[\delta].$ By choosing $\Delta$-invariant
$\O_l$-lattices $T_l$ of $V_l$ and a basis of $T_l$ the above
isomorphism induces an embedding of $\O$-algebras
\begin{eqnarray}
\label{emb}\Omega:\O[\Delta]\subseteq \prod_{l=1}^k M_{n_l}(\O_l),
\end{eqnarray}whose cokernel
is finite and annihilated by the order $|\Delta_p|$ of a $p$-Sylow
group $\delta_p$ of $\Delta.$ I am grateful to D. Vogel for pointing out to me the following

\begin{ex} \label{p-group}
If $\Delta$ is a  $p$-group, by a theorem of Fong one obtains an
embedding of $\zp$-algebras with finite cokernel
\[\zp[\Delta]\subseteq \prod_{l=1}^k M_{n_l}(\zp[\zeta_{p^{m_l}}]),\]
i.e.\ $K_l=\qp(\zeta_{p^{m_l}})$ where $\zeta_{p^{m_l}}$ denotes a
primitive $p^{m_l}$th root of unity for some $m_l.$
\end{ex}

In order to define equivariant Euler-characteristics we will first
discuss  which different possibilities of   (relative) $K$-groups
exist in which they could live. To this end consider the following
 diagram

\[\xymatrix{
  K_0(\F(\O[\Delta]))\ar[ddd]_-{forget} & K_0(\O[\Delta],\qp)\ar@{->>}[r]^-{\theta}\ar[l]_{\psi}\ar[ddd]_{forget} & K_0(\prod_{l=1}^k M_{n_l}(\O_l),\qp) \ar@{=}[d] \\
  &   & {\prod_{l=1}^k K_0(M_{n_l}(\O_l),\qp)\ar@{=}[d]_{Morita} }\\
    &   &    {\prod_{l=1}^k K_0(\O_l,\qp)\cong{\prod_{l=1}^k
  \z} \ar@{=}[d]^{\prod [\kappa_l:\kappa]}_{forget} }\\
  K_0(\F(\O))\ar@{=}[r] & K_0(\O,\qp)\cong\z &   {\prod_{l=1}^k K_0(\O,\qp)\cong \prod_{l=1}^k
  \z\ar[l]_-{\sum'}
}}\]

Here we write $ \F(\O[\Delta])$ and $\F(\O)$ for the categories of
finite $\O[\Delta]$- and $\O$-modules, respectively. For a
$\zp$-algebra $A$ which is free and finitely generated as
\zp-module we denote by $K_0(A,\qp)$ the relative $K$-group
associated to the ring homomorphism $A\to A\otimes_\zp \qp$ in the
sense of Swan \cite{swan}.  But we mention that  $K_0(A,\qp)$ can
be identified with the Grothendieck group of the category of
finite $A$-modules of finite projective dimension. Therefore the
exact functor
\[\{\mbox{finite $A$-modules of finite projective
dimension}\}\to\{\mbox{finite $A$-modules}\}\] induces a
homomorphism $\psi,$ but note that $K_0(\O[\Delta],\qp)$ and
$K_0(\F(\O[\Delta])$ are not isomorphic in general.

Since $\O$ is regular we have the first isomorphism in the bottom
line which is induced by the analogue of $\psi.$ The maps labelled
``forget" are induced by the forgetful functor
$\mbox{$\O[\Delta]$-modules}\to\mbox{$\O$-modules}$ and
$\mbox{$\O_l$-modules}\to\mbox{$\O$-modules},$ respectively. Thus
the first quadrant obviously commutes. In the second quadrant the
map $\theta$ is induced by the base change functor $\prod_{l=1}^k
M_{n_l}(\O_l)\otimes_{\O[\Delta]}-,$ it is surjective because of
the canonical isomorphism  $K_1(K[\Delta])\cong K_1(\prod_{l=1}^k
M_{n_l}(K_l))$ (use the (long) exact sequence of relative
$K$-theory and functoriality of base change). The map labelled
``Morita" comes from Morita equivalence, more precisely it is
induced by the functor $ T_l\otimes_{M_{n_l}(\O_l)} -,$ where
$T_l$ is considered as $\O_l-M_{n_l}(\O_l)$- bimodule in the
obvious way. Note that we have natural isomorphisms of
$\O_l-M_{n_l}(\O_l)$- bimodules
\[M_{n_l}(\O_l)\cong\bigoplus_{i=1}^{n_l} e^i_{n_l}M_{n_l}(\O_l)\cong
T_l^{n_l},\] where $e^i_{n_l},$ $1\leq i\leq n_l,$ denotes the
idempotent
\[e^i_{n_l}\in M_{n_l}(\O_l)\] whose unique non-zero entry is $1$
at the $i$th diagonal place. Thus the functor $
T_l\otimes_{M_{n_l}(\O_l)} -$ is equivalent to applying $e^i_{n_l}
$ for some $i,$ say $1.$

Identifying $K_0(\O_l,\qp)$ with $\z$ by associating to a finite
$\O_l$-module $M$ its \linebreak $\mathrm{length}_{\mathcal{O}_l}(M)$ and
analogously for $K_0(\O,\qp)$ one immediately verifies that the
forgetful functor induces the map \[[\kappa_l:\kappa]:\z\to\z,
a\mapsto [\kappa_l:\kappa]a,\] where $\kappa_l$ and $\kappa$
denote the residue class fields of $\O_l$ and $\O,$ respectively.


Finally, by definition the map $\sum':\prod_{l=1}^k
  \z\to\z$ maps $(a_l)_l$ to $\sum_{l=1}^k {n_l}\cdot a_l $ and we
  claim

\begin{lem}\label{diag}
The above diagram is commutative.
\end{lem}

\begin{proof}
Note that the group $K_0(\O[\Delta],\qp)$ can be generated by triples \linebreak $(\O[\Delta],\O[\Delta],f\otimes_\zp\qp)$ where
$f:\O[\Delta]\to\O[\Delta]$ is a $\O[\Delta]$-module homomorphism
with finite kernel and cokernel. The image of such a class in $K_0(\O,\qp)\cong\z$
is $\mathrm{length}_{\mathcal{O}}(\mathrm{\coker}(f)).$

Using the embedding \ref{emb} one
 obtains the following commutative diagram of $\O$-modules with exact rows

 \[\xymatrix{
   0\ar[r] & {\O[\Delta]\ar[r]\ar[d]^f }& {\prod_{l=1}^k M_{n_l}(\O_l)\ar[r]\ar[d]^{F}} & C\ar[r]\ar[d] & 0\phantom{,} \\
   0\ar[r] & {\O[\Delta] \ar[r]}& {\prod_{l=1}^k M_{n_l}(\O_l)\ar[r]} & C\ar[r] & 0, \
 }\]
  where $F=\prod_{l=1}^k M_{n_l}(\mathcal{O}_l)\otimes_{\mathcal{O}[\Delta]}
 f.$ Since $C$ is finite and observing that $\ker(f)=\ker(F)=0$ one concludes by the snake lemma that
 \[\mathrm{length}_{\mathcal{O}}(\mathrm{\coker}(f))=\mathrm{length}_{\mathcal{O}}(\mathrm{\coker}(F)).\]
On the other hand $F=\sum_{l=1}^k
T_l^{n_l}\otimes_{\mathcal{O}[\Delta]} f$ and thus
\[\mathrm{length}_{\mathcal{O}}(\mathrm{\coker}(F))=\sum_{l=1}^k n_l\cdot\mathrm{length}_{\mathcal{O}}(\mathrm{\coker}(F_l)),\]
where $F_l=T_l\otimes_{\mathcal{O}[\Delta]} f.$ This implies the
lemma.
\end{proof}

Now we are in a position to discuss different
Euler-characteristics: Let $M$ be a finitely generated $\O \kl
G\kr$-module with finite projective dimension (if $G$ does not have an element of order $p$ this holds for every $\O \kl
G\kr$-module) and assume that all homology groups $\H_i(U,M)$ are
finite. With other words, choosing a projective resolution
$P^\bullet\to M$ we obtain a bounded complex of finitely generated
projective $\O[\Delta]$-modules
\[P_U^\bullet=\O[\Delta]\otimes_{\mathcal{O} \kl
G\kr}^{\mathbb{L}}M\] with finite homology groups. Using as before
the category of virtual objects we obtain an element
$(P_U^\bullet,\lambda_{P_U^\bullet})\in V(\O[\Delta], K[\Delta])$
whose class in $K_0(\O[\Delta]),\qp)$  we denote by
\[\chi_u(U,M)\] because it should be considered as the
``universal" $U$-Euler characteristic in this situation. On the
other hand we define \begin{eqnarray}\chi_f(U,M)&:=&\sum_i
(-1)^{i}\;[\H_i(P_U^\bullet)]\in K_0(\F(\O[\Delta]))\\
&=&\sum_i (-1)^{i}\;[\mbox{Tor}_i^{\mathcal{O}\kl G
\kr}(\O[\Delta],M)]\\
&=&\sum_{i} (-1)^{i}\;[\H_i(U,M)].
\end{eqnarray}
If it happens that all $\H_i(P_U^\bullet)$ are in
$D^p(\O[\Delta]),$ then it follows from \eqref{det-homology} that
\[\psi(\chi_u(U,M))=\chi_f(U,M),\] but we don't know whether there is any  
relation  in general. In any case their
images under the forgetful functor coincide by the regularity of
$\O$ and \eqref{det-homology} and they equal the ``absolute"
$U$-Euler characteristic
\begin{eqnarray}
\chi_{a}(U,M)&:=&\sum_i
(-1)^{i}\;\mathrm{length}_\mathcal{O}(\H_i(U,M))\in \z\cong
K_0(\O,\qp).
\end{eqnarray}

Considering the right hand side of the above diagram we are lead
to consider $G$-Euler characteristics of twisted modules. For any finitely generated free $\O$-module with continuous $G$-action we
define \begin{eqnarray} \phantom{mm}\chi_{a}(G,T,M)&:=&\sum_i
(-1)^{i}\;\mathrm{length}_\mathcal{O}(\mbox{Tor}_i^{\mathcal{O}\kl
G \kr}(T,M))\in \z\cong K_0(\O,\qp),
\end{eqnarray}
if all Tor-groups are finite, compare with \cite[{\S}3]{howson2001} where it is shown that this {\em twisted Euler characteristic} only depends on $V=T\otimes_\O K$ if $G$ is a pro-$p$-group, a fact which will turn out automatically in the situation we will consider in section \ref{descent}. For the trivial representation we
also write \[\chi_{a}(G,M):=\chi_{a}(G,\O,M),\] which is conform
with our definition of $\chi_{a}(U,M)$ for $U=G.$

From Lemma \ref{diag} and its proof we obtain

\begin{prop}
The Euler characteristic $\chi_u(U,M)$ is defined if and only if $\chi_{a}(G,T_l,M)$
is defined for all $1\leq l\leq k.$ If these equivalent statements
hold, then the image of $\chi_u(U,M)$ under the above map
$K_0(\O[\Delta],\qp)\to \prod_{l=1}^k K_0(\O,\qp)$ equals
$(\chi_{a}(G,T_l,M))_l.$ Consequently we have
\[\chi_{a}(U,M)=\sum_{l=1}^k n_l\cdot\chi_{a}(G,T_l,M)=\chi_{a}(G,M) + \sum_{l\neq 1}n_l\cdot\chi_{a}(G,T_l,M) ,\]
where we assume without loss of generality $T_1=\O.$
\end{prop}

Let us consider special cases.

{\bf Case I:} $\Delta$ is a $p$-group.

Then $\O[\Delta]$ is a local ring  with unique simple module
$\kappa$ (with trivial action) \cite[ch.\ V]{nsw}. Thus
$K_0(\F(\O[\Delta])=K_0(\F(\O))=\z$ and therefore $\chi_f(U,M)$
and $\chi_{a}(U,M)$ coincide. By example \ref{p-group} we may take
$\O=\zp$ here and then $[\kappa_l:\kappa]=1$ for all $l.$


{\bf Case II:} $\Delta$ is of order prime to $p.$

Now the embedding $\Omega$ is an isomorphism \cite[prop.\
27.1]{cr} and thus $\chi_u(U,M)$ is completely determined by the
tuple $(\chi_{a}(G,T_l,M))_l.$ Moreover, $\O[\Delta]$ is  regular
now and thus $\psi:K_0(\O[\Delta],\qp)\cong K_0(\F(\O[\Delta]))$
is an isomorphism.

\POP
\PUSH{k-theory.tex}%
\Section{Descent of $K$-theory}\label{desccent}

Let $G$ be  a compact $p$-adic Lie group satisfying the following
conditions: (i) $G$ is pro-$p,$ (ii) $G$ has no element of order
$p,$ and (iii) $G$ has a closed normal subgroup $H$ such that
$\Gamma:=G/H\cong\zp.$ In particular, $G$ is isomorphic to the
semi-direct product $H\rtimes \Gamma.$ By \mbox{\rm $\La(G)$-mod} and
$\mbox{\rm $\La(G)$-mod}^H$ we denote as before the category of finitely
generated  $\La(G)$-modules and its full subcategory consisting of
modules which are finitely generated as
$\La(H)$-module, respectively. The latter is a full subcategory of
the category of finitely generated torsion $\La(G)$-modules.

In \cite[{\S}4]{cs2} Coates-Schneider-Sujatha define the alternating
characteristic ideal of  $M\in \mbox{\rm $\La(G)$-mod}^H$  as follows:

\[ \mathrm{Ak}_G(H,M):=\prod_{i\geq 0} \mathrm{char}_\Gamma(\H_i(H,M))^{(-1)^i}.\]
This can be considered as a relative Euler characteristic and following J. Coates we call it the {\em Akashi series} of $M.$ 
This is a fractional ideal of the quotient field $Q(\Gamma)$ of
$\La(\Gamma)$ and can alternatively be interpreted as an element
of $Q(\Gamma)^\times/\La(\Gamma)^\times.$  Then the above invariant
induces the following map of $K$-groups
\[\mathrm{Ak}_G(H,-):K_0(\mbox{\rm $\La(G)$-mod}^H)\to K_0(\Lambda(\Gamma),
Q(\Gamma))\cong Q(\Gamma)^\times/\La(\Gamma)^\times,\] where
$K_0(\Lambda(\Gamma), Q(\Gamma))$ denotes the relative $K$-group
of the ring-homomorphism $\La(\Gamma)\to Q(\Gamma)$ in the sense
of Swan \cite{swan}, see section \ref{swan}.

Now assume that the multiplicative set $\T$ defined in subsection \ref{Ore-sets} is an Ore-set and recall from section \ref{swan} that we had associated with every $M\in \La(G)\mbox{-mod}^H=\La(G)\mbox{-mod}_{\T-\mathrm{tor}}$ an characteristic class \[\mathrm{char}_G(M):=\mathrm{char}_{\La(G)}(M)\] in the group $K_0(\La,\La_\mathcal{T})$ which 
 fits into the following short exact sequence
\[ \xymatrix{
  K_1(\La) \ar[r]\ar@{=}[d] & K_1(\La_\mathcal{T})\ar[r]\ar@{=}[d] & K_0(\La,\La_\mathcal{T}) \ar[r]\ar@{=}[d]& 0 \\
  {\La^\times/[\La^\times,\La^\times]\ar[r]} & {\La_\mathcal{T}^\times/[\La_\mathcal{T}^\times,\La_\mathcal{T}^\times]\ar[r]} & {\La_\mathcal{T}^\times/[\La_\mathcal{T}^\times,\La_\mathcal{T}^\times]\La^\times\ar[r] }&
  1.} \]
Thus we  considered $\mathrm{char}_G(M)$ also as an element $F_M
[\La_\mathcal{T}^\times,\La_\mathcal{T}^\times]\La^\times\in
\La_\mathcal{T}^\times/[\La_\mathcal{T}^\times,\La_\mathcal{T}^\times]\La^\times$
and we called (any such choice) $F_M\in\La_\mathcal{T}^\times$ a {\em characteristic}
element of $M.$

Now  consider the canonical ring homomorphism
$\pi_H:\La(G)\to\La(\Gamma)$ which is induced by the group
homomorphism $G\to\Gamma.$ It 
 induces a commutative diagram
of $K$-groups with exact rows

\[  \xymatrix{
  K_1(\La) \ar[r]\ar@{->}^{(\pi_H)_\ast}[d] & K_1(\La_\mathcal{T})\ar[r]\ar@{->}^{(\pi_H)_\ast}[d] & K_0(\La,\La_\mathcal{T}) \ar[r]\ar@{->}^{(\pi_H)_\ast}[d]& 0 \\
  K_1(\La(\Gamma)) \ar[r]\ar@{=}[d] & K_1(Q(\Gamma))\ar[r]\ar@{=}[d] & K_0(\La(\Gamma),Q(\Gamma)) \ar[r]\ar@{=}[d]& 0 \\
  {\La(\Gamma)^\times \ar[r] }& Q(\Gamma)^\times\ar[r] & Q(\Gamma)^\times/\La(\Gamma)^\times\ar[r] &
  1,
}\]

where the middle map $(\pi_H)_\ast$ is induced by
$\La_\mathcal{T}^\times\to Q(\Gamma)^\times.$ The following result is
now almost self-proving.

\begin{prop}\label{descent}
Let $M\in\mbox{\rm $\La(G)$-mod}^H.$ Then  the following
holds:\[\mathrm{Ak}_G(H,M)\equiv\pi_H(\mathrm{char}_G(M))\equiv\pi_H (F_M)
 \mbox{ \em mod } \La(\Gamma)^\times,\] i.e.\
there is a commutative diagram
\[\xymatrix{
  {K_0(\mbox{\rm $\La(G)$-mod}^H) \ar@^{=}[r]\ar_{\mathrm{Ak}_G(H,-)}[dr]}&   {K_0(\La(G),\La(G)_\mathcal{T}) \ar^{(\pi_H)_\ast}[d]} \ar@^{=}[r]&  {\La(G)_\mathcal{T}^\times/[\La(G)_\mathcal{T}^\times,\La(G)_\mathcal{T}^\times]\La(G)^\times\ar^{(\pi_H)_\ast}[dl]}.\\
 {}  &  {K_0(\La(\Gamma),Q(\Gamma))} & {} }\]
\end{prop}

\begin{proof}
The second $\equiv$ follows from the functoriality of our
identifications. To prove the first choose a projective
resolution $P^{\bullet} \to M$ by a perfect complex $P^\bullet.$
As explained in the paragraph before Proposition \ref{ident} $\mathrm{char}_G(M)$ is given by the class of $([P^\bullet],\lambda_{P^\bullet})\in V(\La,\La_\T)$ where \[\lambda_{P^\bullet}:\La_\T\otimes_\La [P^\bullet]=[\La_\T\otimes_\La P^\bullet]\cong \mathbf{1}_{V(\La_\mathcal{T})}\] is the canonical isomorphism of virtual objects in $V(\La_\mathcal{T})$ associated to the quasi-isomorphism $\La_\T\otimes_\La P^\bullet\to 0.$ Now the base chance maps $\pi_H:\La\to \La(\Gamma)$ and $\pi_H:\La_\T\to Q(\Gamma)$ induce a morphism of Picard categories \[(\pi_H)_\ast:V(\La,\La_\T)\to V(\La(\Gamma),Q(\Gamma))\] under which $([P^\bullet],\lambda_{P^\bullet})$ is mapped to $([\La(\Gamma)\otimes_\La P^\bullet], Q(\Gamma)\otimes_{\La_\T} \lambda_{P^\bullet})$ where we use the compatibility of (non-commutative) determinants $[-]$ with arbitrary change of rings: \[Q(\Gamma)\otimes_{\La_\T}\La_\T\otimes_\La [P^\bullet]=Q(\Gamma)\otimes_{\La(\Gamma)}[\La(\Gamma)\otimes_{\La}P^\bullet],\]
\[Q(\Gamma)\otimes_{\La_\T} \mathbf{1}_{V(\La_\mathcal{T})}=\mathbf{1}_{V(Q(\Gamma))}\]
and thus 
\[\xymatrix{{Q(\Gamma)\otimes_{\La_\T} \lambda_{P^\bullet} :}\;
  Q(\Gamma)\otimes_{\La(\Gamma)}[\La(\Gamma)\otimes_{\La}P^\bullet]\ar[r]^-{\cong} &  {\mathbf{1}_{V(Q(\Gamma))}.}
}\]

Since for the commutative rings $\La(\Gamma)$ and $Q(\Gamma)$ we
can replace $[-]$ by the functor $\det(-),$ the above shows that
$\pi_H(\mathrm{char}_G(M))$ in $K_0(\La(\Gamma),Q(\Gamma))$ is
represented by the pair
$(\det_{\La(\Gamma)}(\La(\Gamma)\otimes_{\La} P^\bullet),
\lambda_{\La(\Gamma)\otimes_{\La} P^\bullet}),$ where
\[\lambda_{\La(\Gamma)\otimes_{\La(G)}
P^\bullet}=Q(\Gamma)\otimes_{\La_\T} \lambda_{P^\bullet}:
\mathrm{det}_{Q(\Gamma)}(Q(\Gamma)\otimes_{\La} P^\bullet)\cong
\mathrm{det}_{Q(\Gamma)}(0)= (Q(\Gamma),0).\] Since $\La(\Gamma)$
is regular we obtain by \cite{km} or \cite[(9)]{bf} a canonical
isomorphism
\begin{eqnarray*}
\mathrm{det}_{\La(\Gamma)}(\La(\Gamma)\otimes_{\La(G)}
P^\bullet)&\cong&\bigotimes_{i\in\mathbb{Z}}
\mathrm{det}_{\La(\Gamma)}(H^i(\La(\Gamma)\otimes_{\La(G)} P^\bullet)^{(-1)^{i-1}}\\
&=&\bigotimes_{i\in\mathbb{Z}}
\mathrm{det}_{\La(\Gamma)}(H_i(H,M))^{(-1)^{i-1}}.
\end{eqnarray*}
 But Kato  \cite[prop.\ 6.1]{kato} observed that $\mathrm{char}_{\Gamma}(N)\mathrm{``}=\mathrm{"}(\det_{\La(\Gamma)}(N))^{-1}$ for every torsion $\La(\Gamma)$-module $N,$ in the sense that under the canonical isomorphism \[\mathrm{det}_{Q(\Gamma)}(Q(\Gamma)\otimes_{\La(\Gamma)}N)\cong \mathrm{det}_{Q(\Gamma)}(0)=(Q(\Gamma),0)\] $\mathrm{det}_{\La(\Gamma)}(N)\subseteq \mathrm{det}_{Q(\Gamma)}(Q(\Gamma)\otimes_{\La(\Gamma)}N)$ is mapped to the fractional ideal generated by $\mathrm{char}_{\Gamma}(N)^{-1}.$ This implies that $\det_{\La(\Gamma)}(\La(\Gamma)\otimes_{\La(G)}
P^\bullet)$ is mapped under $\lambda_{\La(\Gamma)\otimes_{\La(G)}
P^\bullet}$ to $\La(\Gamma)\mathrm{Ak}_G(H,M)\subseteq
Q(\Gamma)$ and the proposition follows.

In order to understand better the different constructions we now sketch a second proof which avoids the use of virtual objects (which have nice functorial and universal properties as we have seen above but which are quite difficult to imagine). To this end we represent $M$  by the class $[P^+,P^-,\phi]$ where $\phi:P^+_\T\to P^-_\T$ is an isomorphism obtained by choosing successively sections, see section \ref{swan}. Since $M$ is torsion $P^+$ and $P^-$ are free of the same rank, $r$ say, and thus the class of $\phi\in \mathrm{Aut}_{\La_\T}(\La_\T^r)$ in $K_1(\La_\T)$ is a characteristic element of $M.$ Base change with respect to $\pi_H$ gives a class $[P^+_H,P^-_H,\phi_H]$ where \[\phi_H=Q(\Gamma)\otimes_{\La_\T}\phi:Q(\Gamma)\otimes_{\La(\Gamma)}P^+_H \cong Q(\Gamma)\otimes_{\La(\Gamma)}P^-_H\] is the induced isomorphism of free $Q(\Gamma)$-modules. Of course, $\phi_H=(\pi_H)_\ast(\phi)$ is a characteristic element of the complex $P^\bullet_H.$ Now using Proposition \ref{ident} one can proceed as above using determinants or - if one even wants to avoid them - it is not difficult but tedious to show directly that $[P^+_H,P^-_H,\phi_H]$ can be expressed as alternating sum of classes associated to the homology groups of $P^\bullet_H$ (cf.\  \cite{burns}), which are $\La(\Gamma)$-torsion modules and thus represented by classes \linebreak $[\La(\Gamma),\La(\Gamma),\mathrm{char}_\Gamma(\H^i(H,M))].$ 

But the most elegant way to prove the proposition is the following: Use the fact that $\La(G)_\T^\times$ is generated by the elements in $\T$ and check the commutativity just for $t\in \T.$ Obviously, $t$ is mapped to the module $\La/\La t$ whose Akashi series is just $\pi_H(t).$ 
\end{proof}


\begin{rem}
(The commutative case) Let $G$ be isomorphic to $\zp^d$ for some integer $d\geq 2$ and fix for the moment a subgroup $H$ such that $G/H\cong\zp.$ Recall from \cite[{\S}4.5]{bourbaki} that in this situation we have  canonical isomorphisms \[K_0(\mbox{\rm $\La(G)$-mod}_\mathrm{tor}/\mathcal{ PN})\cong K_0(\mbox{\rm $\La(G)$-mod}_\mathrm{tor})\cong \mathrm{Div}(\La(G)),\]
where $\mathrm{Div}(\La(G))\cong \bigoplus_{P\in \mathcal{P}} \z$ denotes the group of divisors of $\La(G)$ and $\mathcal{P}$ denotes a system of representatives of classes of irreducible elements of $\La(G).$ In particular, using this identification an element $f=(unit)\prod_{\mathcal{P}}P^{v_P(f)} \in Q(G)^\times\cong K_1(Q(G))$ is mapped to its divisor $\mathrm{div}(f)=\sum_{\mathcal{P}}v_P(f)P $ under the connecting map in the localization sequence of $K$-theory. Now one sees immediately that the injection $\La_\T^\times\to Q(G)^\times$ induces a commutative diagram \[\xymatrix{
     {\La_\T^\times/\La(G)^\times}\ar@{^{(}->}[d]\ar@{=}[r] & K_0(\mbox{\rm $\La(G)$-mod}^H) \ar@{^{(}->}[d]\\
  Q(G)^\times/\La(G)^\times \ar@{=}[r]& K_0(\mbox{\rm $\La(G)$-mod}_\mathrm{tor}).
}\] 
In particular, we see that our characteristic class or element for modules which are finitely generated over $\La(H)$  can be identified with the usual one of the full torsion category. In the next remark we will see that this fails in the non-commutative situation because the big commutator subgroups of the units destroy the injectivity. The image of $K_0(\mbox{\rm $\La(G)$-mod}^H) $ in $\mathrm{Div}(\La(G))$ is precisely \[\mathrm{Div}(\La(G))_H:=\bigoplus_{P\in \mathcal{P}_{\mathrm{red}}} \z,\] where $\mathcal{P}_{\mathrm{red}}$ denotes the subset of $\mathcal{P}$ consisting of all elements having finite reduced order (with respect to $H$). In geometric terms $\mathrm{Div}(\La(G))_H$ consists precisely of those cycles, which have codimension $1$ and which have good intersection with the closed subscheme defined by $\m(H).$\\
If now $H$ ranges over all subgroups of $G$ having a quotient isomorphic to $\zp$ one concludes using \cite[lem.\ 2]{green} that \[\dirlim H K_0(\mbox{\rm $\La(G)$-mod}^H))\cong \mathrm{Div}(\La(G)\otimes_\zp \qp)\cong \bigoplus_{p\neq P\in \mathcal{P}} \z,\] i.e.\ up to $\zp$-torsion modules one obtains all $\La(G)$-torsion modules this way.
\end{rem}

\begin{rem} In Remark \ref{ex-css2} we have discussed why we have to replace 
 the full torsion
category $\mbox{\rm $\La(G)$-mod}_\mathrm{tor}$ by a smaller subcategory like
$\mbox{\rm $\La(G)$-mod}^H\cong\linebreak\mbox{$\La$-mod}_{\T\mathrm{-tor}}.$
Note that  the example mentioned there also shows that the canonical map
\[K_0(\mbox{\rm $\La(G)$-mod}^H)\to K_0(\mbox{\rm $\La(G)$-mod})\]
is not injective  in general. In contrast to our expectation from
the commutative theory (e.g.\ Gersten's conjecture) it is also
remarkable that the canonical map
\[K_0(\mbox{$\La$-mod}_{\La(H)\mathrm{-tor}})\to K_0(\mbox{\rm $\La(G)$-mod}^H)\]
is {\em not} trivial, i.e.\ also pseudo-null $\La(G)$-modules (in
$\mbox{\rm $\La(G)$-mod}^H$) can give rise to non-trivial classes.\\
\end{rem}

\begin{rem}
Let $F_M=gh^{-1}$ be a characteristic element for $M\in \La(G)\mbox{-mod}^H$ with $g,h\in \La_\T.$ Then we have \[\rk_{\La(H)} M=\mathrm{ord}^\mathrm{red}(F_M):=\mathrm{ord}^\mathrm{red}(g)-\mathrm{ord}^\mathrm{red}(h),\] where for $g\in\La(G)$ we denote by $\mathrm{ord}^\mathrm{red}(g)$ the order of the power series $\psi_H(g)\in\kappa\kl \Gamma\kr\cong\kappa[[T]]$ (cf.\ \cite{ven-weier}).
In particular, $\mathrm{ord}^\mathrm{red}(F_M)\geq 0$ is independent of the chosen fraction. We think of  $\rk_{\La(H)} M$ as a generalized $\lambda$-invariant. The claim follows from the fact that the $\La(H)$-rank is additive on short exact sequences and thus induces a group homomorphism $\rk_{\La(H)}:K_0(\La(G)\mbox{-mod}^H)\to\z.$ Now $[M]=[\La/\La g]-[\La/\La h]$ and from the Weierstrass preparation theorem (loc.\ cit.) we conclude that $\rk_{\La(H)} \La/\La g=\mathrm{ord}^\mathrm{red}(g).$
\end{rem}

\begin{prop}\label{chardet}
  Assume that $M$ in $\mbox{\rm $\La(G)$-mod}^H$ satisfies a relation
$[M]=[\coker(f)]-[\coker(g)]$ for some $f,g\in M_n(\La(G))$ such
that $\coker(f),\coker(g)$ belong to $\La(G)\mbox{-mod}^H.$
Then for all
  representations $\rho:G\to \mathrm{Aut}_\mathcal{O}(T)$ we have
  \[\mathrm{Ak}_G(H,T\otimes_\O M)=\mathrm{det}_{Q(\Gamma)}(\pi_H(\mathrm{tw}_\rho(f))\cdot \pi_H(\mathrm{tw}_\rho(g))^{-1})\]
  in $K_0(\O\kl \Gamma\kr,Q(\Gamma)).$ In particular, it holds for any choice of $F_M$ that \[\mathrm{Ak}_G(H,T\otimes_\O M)=\mathrm{det}_{Q(\Gamma)}(\pi_H(\mathrm{tw}_\rho(F_M)).\]
  \end{prop}

\begin{proof}
Applying $\rho_\ast$ to the relation
$[M]=[\coker(f)]-[\coker(g)]$ gives the relation
\[[T\otimes_\O
M]=[\coker(\mathrm{tw}_\rho(f))]-[\coker(\mathrm{tw}_\rho(g))]\]
which in turn induces after evaluating the functor
$\mathrm{Ak}_G(H,-)$ (with $\O$-coefficients)
\begin{eqnarray*}
\mathrm{Ak}_G(H,T\otimes_\O
M)&=&\mathrm{Ak}_G(H,\coker(\mathrm{tw}_\rho(f)))\cdot \mathrm{Ak}_G(H,\coker(\mathrm{tw}_\rho(f)))^{-1}\\
&=&\mathrm{det}_{Q(\Gamma)}(\pi_H(\mathrm{tw}_\rho(f)))\cdot\mathrm{det}_{Q(\Gamma)}(
\pi_H(\mathrm{tw}_\rho(g)))^{-1}\\
&=&\mathrm{det}_{Q(\Gamma)}(\pi_H(\mathrm{tw}_\rho(f))\cdot
\pi_H(\mathrm{tw}_\rho(g))^{-1}),
\end{eqnarray*}
where the second equality is easily verified, compare with the
proof of Lemma \ref{twisting-lem}.
\end{proof}

In subsection \ref{eq-EC} we have defined an additive Euler characteristic, but for our arithmetic applications the following multiplicative version is more suitable \[\chi(G,T,M):=(\#\kappa)^{\chi_a(G,T,M)} \mbox{ (if } \chi_a(G,T,M)\mbox{ is defined).}\] Similarly we write $\chi(U,M)$ for the multiplicative version of $\chi_a(U,M).$

\begin{prop}\label{EC-evaluation} Assume that $M$ in $\mbox{\rm $\La(G)$-mod}^H$ satisfies a relation
$[M]=[\coker(f)]-[\coker(g)]$ for some $f,g\in M_n(\La(G))$ such
that $\coker(f),\coker(g)$ are finitely generated over $\La(H).$
Then, if $\chi(G,T,M)$ is finite,  $fg^{-1}(\rho)$ is defined (and non-zero)
and we have
\[\chi(G,T,M)=|fg^{-1}(\rho)|_p^{-[K:\qp]},\] where the norm is normalized by $|p|_p=\frac{1}{p}.$
In particular, we have for any choice of $F_M$ 
\[\chi(G,T,M)=|F_M(\rho)|_p^{-[K:\qp]}.\]
\end{prop}

This proposition in junction with Lemma \ref{lattice-indep} shows that $\chi(G,T,M)$ actually depends only on the $p$-adic representation $V=T\otimes_\O K$ and thus we will sometimes write $\chi(G,V,M)$ for it.

\begin{proof}
Under the assumptions $\mathrm{Ak}_G(H,T\otimes_\O M)(0)$ is
defined in the sense of \cite[{\S}4]{cs2} and $\chi(G,T\otimes_\O
M)=|\mathrm{Ak}_G(H,T\otimes_\O M)(0)|_p^{-[K:\qp]}$ by
\cite[lem.\ 4.2]{cs2}. Actually, the lemma is only stated for
$\zp$-coefficients but it obviously generalizes to  the case of
coefficients in $\O$ using the well known fact that \[\#(\O/\O f)=|f|_p^{[K:\qp]}\] for every $f\in \O\setminus\{0\}.$ The idea of that lemma is the following: evaluating an element of $K_0(\La(\Gamma,Q(\Gamma))$ at zero means calculating the $\Gamma$-Euler characteristic of that element (if defined), i.e.\ in this case the $\Gamma$-Euler characteristic of $\prod_{i\geq 0} \mathrm{char}_\Gamma(\H_i(H,M))^{(-1)^i},$ which is related to the $G$-Euler characteristic of $M$ via the (almost degenerating) spectral sequence
\[\H_p(\Gamma,\H_q(H,M))\Rightarrow \H_{p+q}(G,M).\] Inspection of the latter gives their lemma.

By Prop.\ \ref{chardet} we have
$\mathrm{det}_{Q(\Gamma)}(\pi_H(\mathrm{tw}_\rho(f))\cdot
\pi_H(\mathrm{tw}_\rho(g))^{-1})(0)=\mathrm{Ak}_G(H,T\otimes_\O
M)(0)$
  and
thus $fg^{-1}(\rho)$   is defined and the claim follows from the
previous Lemma.
\end{proof}

For $G\cong\zp$ it is well known that also a converse statement holds: If $f_M(0)\neq0,$ then $\chi(G,M)$ is defined. But in higher dimensions one easily constructs counterexamples in which lemma 4.2 of \cite{cs2} and thus this sort of statement fails to be true.

\begin{cor}
Let $U$ be an open normal subgroup of $G$ and $M\in\mbox{\rm $\La(G)$-mod}^H.$ Then, for any choice of $F_M,$ we have
\[\chi(U,M)=\prod_{\rho} |F_M(\rho)|_p^{-[K:\qp]n_\rho}\]
where $\rho$ runs through a system of representatives of irreducible $K$-representations of $\Delta:=G/U$ and where $K$ is a splitting field of $\Delta.$
\end{cor}

This result should be compared with  well-known formulas in the $\zp$-situation, see \cite[Ch.\ V {\S} 3 Ex.\ 3]{nsw}.

At the end of this section we want to determine characteristic elements in several examples and use them to calculate the Euler characteristics. Of particular interest is  D. Vogel's example of a non-principal reflexive left ideal in a Iwasawa algebra which we state first:

\begin{ex}\label{vogel}
(\cite[Appendix A.3]{ven-weier}) Consider the Iwasawa algebra $\Lambda=\zp\kl G\kr$ of the semidirect product
$G=\Gamma_1\rtimes\Gamma_2$ of two copies $\Gamma_1,\Gamma_2$ of
$\zp$ where $p$ is any odd prime number.  Let the action of $\Gamma_2$ on $\Gamma_1$
be given by the continuous group homomorphism
$$
\rho:\Gamma_2\rightarrow{\rm Aut}(\Gamma_1), \ \gamma_2\mapsto(\gamma_1\mapsto\gamma_1^{p+1})
$$
for appropriate generators $\gamma_i$ of $\Gamma_i$.
As in \cite{ven-weier} we identify $\La$ with the skew power series ring
$$
\zp\kl G\kr\cong \zp[[X,Y;\sigma,\delta]]
$$
where $X:=\gamma_1-1,Y:=\gamma_2-1$ and recall that the ring automorphism $\sigma$ of $\zp[[X]]$ is induced
by $X\mapsto(X+1)^{p+1}-1$ and $\delta$ is the $\sigma$-derivation given by $\delta=\sigma-\mbox{id}$.
Now let $u$ be in $\zp$ with the following properties: 
(i) $u\pi+\sigma^2(\pi)$ is divisible by $\sigma(\pi)$ in $R$, (ii) 
  $u\equiv 1\bmod p$ and (iii)  $\frac{u\pi+\sigma^2(\pi)}{\sigma(\pi)}\equiv 2\bmod (p,X).$ 
 Putting $\xi=X-p$ he defines skew polynomials \[f:=Y^2 + (2-\frac{u\xi+\sigma^2(\xi)}{\sigma(\xi)})Y +(u-\frac{u\xi+\sigma^2(\xi)}{\sigma(\xi)}+1)\] and \[k:=\xi Y+(\xi-\sigma(\xi)).\] He then shows that an element $u\in\zp$ with the above properties exists  and that the left ideal $L:=\La f+\La k $ generated by $f$ and $k$ is reflexive, but  cannot be generated by a single element in $\La.$ \\
Concerning this example   D. Vogel  shows in \cite{vo} that there is an exact sequence of $\La$-modules \[\xymatrix@1{
   {\ 0 \ar[r] } &  {\ \La/\La N \ar[r]^{d} } &  {\ \La/\La f \ar[r] } &  {\ \La/L \ar[r] } &  {\ 0, } 
}\]
where $N=Y+1-u$ and the map $d$ is given by $\lambda +\Lambda N\mapsto \lambda k +\Lambda f$ for $\lambda\in \La.$ Thus setting \[F_{\La/L}=fN^{-1}=\{Y^2 + (2-\beta)Y +(u-\beta+1)\}\{Y+1-u\}^{-1}\] where \[\beta=\frac{u\xi+\sigma^2(\xi)}{\sigma(\xi)}\] is a characteristic element for $\La/L$ and $ \mathrm{Ak}_G(H,\La/L)$ is generated by  \[\pi_H(F_{\La/L})=\frac{Y^2+(1-u)Y}{Y+1-u}=Y.\] In particular, $\chi(G,\La/L)$ is not defined. \end{ex}

\begin{ex}\label{ex-char}
\begin{enumerate}
\item As in   Example \ref{vogel} let $G\cong\zp\rtimes\zp,$ where now the action is defined by an arbitrary non-trivial continuous  character $\rho :\zp\to \mathbb{Z}_p^\times.$ In particular, we have $\gamma h\gamma^{-1}=h^\delta$ for fixed topological generators $h,\gamma$ and with non-trivial $\delta=\rho(\gamma)\in\mathbb{Z}_p^\times.$ Again we identify $\La=\La(G)$ with $\zp[[X,Y;\sigma,\delta]],$ where $X=h-1$ and $Y=\gamma-1.$ In this case it is a little exercise to deduce that the pro-$p$ Fox-Lyndon resolution (\cite[5.6.6]{nsw}, see also \cite[Satz 7.7]{koch}) $P^\bullet \to \zp$ has the following form \[\xymatrix@1{
   {\ 0 \ar[r] } &  {\ \La \ar[r]^{d_2} } &  {\ \La^2 \ar[r]^{d_1} } &  {\ \La \ar[r] } &  {\ \zp \ar[r] } & {\ 0}, 
}\]
where $d_2$ is given by right multiplication with the matrix $\left(\begin{array}{cc}
M & N 
\end{array}\right)$ where 
$N=-\sigma(X)$ and $M=Y+ \frac{X-\sigma(X)}{X}$ (or equivalently, if $\delta$ is a positive integer, $N=1-\gamma h\gamma^{-1}$ and $M=\gamma -\sum_{i=0}^{\delta-1}h^i$). The map $d_1$ corresponds to the matrix $\left(\begin{array}{c}
  X \\
  Y \\
\end{array}\right).$ Tensoring with $\La_\T$ gives a short split exact sequence \[\xymatrix@1{
  P^\bullet_\T: & {\ 0 \ar[r] } &  {\ \La_\T \ar[r]^{d_2}  } &  {\ \La_\T^2 \ar[r]^{d_1}\ar@/^/[l]^{t} } &  {\ \La_\T \ar[r]\ar@/^/[l]^{s} } &  {\ 0 } .
}\]
A possible choice for $s$ is for instance given by the matrix $(0\;Y^{-1})$ and then $d\circ t=\mathrm{id}-s\circ d_1$ is given by \[\left(\begin{array}{cc}
  1 & -XY^{-1} \\
  0 & 0 \\
\end{array}\right).\] We conclude that thus $t$ corresponds to $\left(\begin{array}{c}
  M^{-1} \\
  0 \\
\end{array}\right).$ Hence, the map \[\xymatrix{
  P^{-}_\T\cong\La_\T^2\ar[r]^{d_1\oplus t} & {\La_\T^2\cong P^{+}_\T}
}\] is given by the matrix \[\left(\begin{array}{cc}
  M^{-1} & X \\
  0 & Y \\
\end{array}\right)\] and using  Proposition \ref{ident}  we see that \[F_\zp=M^{-1}Y=\Big(Y+ \frac{X-\sigma(X)}{X}\Big)^{-1}Y\] is a characteristic element for $\zp.$ For any continuous character $\psi:\zp\to \mathbb{Z}_p^\times$ we obtain characteristic elements \begin{eqnarray*}
F_{\zp(\psi)}&=&\mathrm{tw}_\psi(F_\zp)\\
&=&\Big(\mathrm{tw}_\psi(Y)+\frac{X-\sigma(X)}{X}\Big)^{-1}\mathrm{tw}_\psi(Y)
\end{eqnarray*}
where \[\mathrm{tw}_\psi(Y)=\psi(\gamma)^{-1}(Y+1)-1\]
(viewing $\psi$ via $G\twoheadrightarrow \zp$ as a character of $G$). Thus $ \mathrm{Ak}_G(H,\zp(\psi) )$ is generated by \[\pi_H(F_{\zp(\psi)})=\frac{\psi(\gamma)^{-1}(Y+1)-1}{\psi(\gamma)^{-1}(Y+1)-\rho(\gamma)+1}\] and hence we have \[\chi(G,\zp(\psi))=|\frac{1-\psi(\gamma)}{1-\psi(\gamma)(\rho(\gamma)-1)}|_p=|\psi(\gamma)-1|_p,\] which is different from $1$ for non-trivial $\psi.$
\item Finally we complete the discussion of Example \ref{ex-css}: Using  sequence \eqref{def-seq} to calculate the $H$-homology Coates-Schneider-Sujatha show that the ideal $\mathrm{Ak}_G(H,M)$ is generated by \[f_M(T)=\frac{T-\omega(0)+1}{T-u(0)+1},\] where for any $z\in \La(H)$ we write $z(0)$ for the image of $z$ under the augmentation map in $\zp.$
Since $\phi$ is injective we have \[\omega(0)=1+p^r\neq \phi(h_2)+p^r=u(0)\] and thus $f_M$ is not a unit in $\La(\Gamma).$  From the short exact sequence \eqref{def-seq} one sees immediately that \[F_M=(c-\omega)(c-u)^{-1}\in \La_\T\] is a characteristic element of $M.$ This gives a second calculation of $f_M$ by \ref{descent}. Finally, by \ref{EC-evaluation} the $G$-Euler characteristic is given by \[\chi(G,M)=|\frac{\phi(h_2)-1+p^r}{p^r}|_p,\] which is generically  non-trivial. 
\end{enumerate}
\end{ex}
\POP
\PUSH{charsel.tex}%
\Section{Characteristic elements of Selmer groups}\label{selmer}

In this section we are going to apply the techniques developed so far to study properties of the Selmer group of an elliptic curve over a $p$-adic Lie extension $k_\infty.$ Needless to say that we could also take arbitrary abelian varieties or motives instead, but all the phenomena we want to discuss occur already for elliptic curves, for which moreover lots of examples have been discussed recently (e.g.\ \cite{coates-howsonII}, \cite{cs2}, \cite{hachi-ven},\cite{howson2000}, \cite{howson2001}, \cite{ochi-ven}).
We shall consider the characteristic element associated with the Pontryagin dual of the Selmer group over $k_\infty$ (provided this module is finitely generated over $\La(H),$ where $H$ denotes the Galois group $G(k_\infty/k_{cyc})$) and discuss its relation with the characteristic polynomial of $E$ over the cyclotomic $\zp$-extension $k_{cyc}.$ To this end we first have to recall some facts from the latter theory.

We fix an odd prime $p.$ Let $k$ be a number field, $S$ a finite
set of places of $k$ containing the set $S_p$ of places lying
above $p$ and the set $S_\infty$ of infinite places. By $k_S$ we
denote the maximal outside $S$ unramified extension of $k$ and,
for any intermediate extension $k_S|L|k,$ we write
$G_S(L):=G(k_S/L)$ for the Galois group of $k_S$ over $L.$
 Suppose that
$G_S(k)$ acts continuously and linearly on a vector space $V$ over
$\Qp$ of dimension $d.$ Let $T$ be a Galois invariant
$\zp$-lattice in $V.$ Then $A=V/T$ is a discrete $G_S(k)$-module
which is isomorphic to $(\Qp/\zp)^d$ as an $\zp$-module. We set $V^\ast=\Hom(V, \qp(1)),$ $T^\ast=\Hom(T,\zp(1)))$ and
$A^\ast=\Hom(T,\mu_{p^\infty}).$ Then it is easy to see that
$A^\ast\cong V^\ast/T^\ast.$

Let $k_{cyc}$ denote the cyclotomic $\zp$-extension of $k$ and set
$\Gamma:=G(k_{cyc}/k).$ For an arbitrary number field $L$ and a place $\nu$ of $L$ we
denote by $L_\nu$ the completion of $L$ at $\nu.$ If $L$ is an
infinite extension of $\Q$ we write $L_\nu$ for the limit of the
completions of the finite subextensions of $L$ with respect to the
induced valuations. Note that the decomposition groups
$\Gamma_\nu:=G(k_{cyc,\nu}/k_\nu)$ have finite index in $\Gamma.$

\Subsection{Local Euler factors}\label{local-euler}

In the context of the Selmer group of $A$ over $k_{cyc}$ the
following local cohomology groups show up:
$\H^1(k_{cyc,v},A)^\vee$ and its global version
$\Ind_\Gamma^{\Gamma_\nu} \H^1(k_{cyc,v},A)^\vee.$ They are
finitely generated $\La(\Gamma_\nu)$- and $\La(\Gamma)$-modules,
respectively. {\em For the rest of this section we shall assume
that $\nu$ is {\em not} lying above $p.$} Then we will see that
the above modules are torsion and the  aim of this subsection is to
determine their characteristic ideals. I am very grateful to Yoshitaka Hachimori for discussions on this problem and for pointing out to me that 
Greenberg and Vatsal \cite{greenvat} have given a description of the characteristic power series in a similar context. We follow closely their approach:

The structure of $\H^1(k_{cyc,v},A)^\vee$ is studied in \cite[Prop
2]{green-adic} and also in \cite{ochi-ven2}, \cite{ven}. In
particular, it is known that their $\mu$-invariant is zero and
that there is a canonical Galois equivariant isomorphism
\[\H^1(k_{cyc,\nu},A)^\vee\otimes_\zp \Qp\cong V^\ast(k_{cyc,\nu})\subseteq (V^\ast)^{I_\nu},\]
where $I_\nu$ denotes the inertia subgroup of $G_{k_\nu}.$

By $\mathrm{Frob}_\nu\in G(k_\nu^{nr}/k)$ we denote the
(arithmetic) Frobenius automorphism where $k_\nu^{nr}$ denotes the
maximal unramified extension of $k_\nu,$ which contains
$k_{cyc,\nu}$ because $\nu \nmid p.$ Furthermore, we write
$\gamma_\nu$ for the image of $\mathrm{Frob}_\nu$ in $\Gamma_\nu.$

Let $\alpha_1,\ldots,\alpha_{e_\nu}$ denote the eigenvalues of
$\mathrm{Frob}_\nu$ (counting with multiplicities) acting on
  on the maximal quotient $V_{I_\nu}$ of $V$ on which $I_\nu$
acts trivially; i.e.\ $e_\nu=\dim_\qp V_{I_\nu}.$ Then the
eigenvalues of $\mathrm{Frob}_\nu$ acting on
\[(V^\ast)^{I_\nu}\cong\Hom_\qp(V_{I_\nu},\qp(1))\]
are $q_\nu \alpha_1^{-1},\ldots,q_\nu \alpha_{e_\nu}^{-1},$ where
$q_\nu$ denotes the order of the residue class field of $\nu.$ We
shall apply the next lemma to $W=(V^\ast)^{I_\nu}$ and
$F=\mathrm{Frob}_\nu.$

\begin{lem}
Let $<F>=\widehat{\z}\to GL(W)$ a continuous representation of the
free profinite group with topological generator $F$ on a finite
dimensional $\qp$-vectorspace. Let $F=F_pF_{p'}$ the unique
decomposition of $F$ corresponding to
$\widehat{\z}\cong\zp\times\widehat{\z}_{(p')},$ where $
\widehat{\z}_{(p')}=\prod_{l\neq p} \mathbb{Z}_l.$ Then the
eigenvalues of $F_p$ (counting with multiplicities) acting on
$W^{<F_{p'}>}$ are precisely those eigenvalues of $F$ (counting
with multiplicities) acting on $W$ which are principal units (in
some extension of $\qp$).
\end{lem}

For the proof just note that the image of $<F_{p'}>$ in $GL(W)$ is
a finite group of order prime to $p$ such that $W$ decomposes into
the eigenspaces  of $F_{p'}.$ Of course, $W^{<F_{p'}>}$ is nothing
else than the eigenspace with eigenvalue $1,$ i.e.\ the
eigenvalues of $F$ and $F_p$ coincide on this subspace while on
the other eigenspaces the eigenvalues of $F$ have a non-trivial
prime to $p$ part.

Let \[P_\nu(T)=\det(1-\mathrm{Frob}_\nu T|V_{I_\nu})=\det(1-\mathrm{Frob}_\nu^{-1}q_\nu
T|(V^\ast)^{I_\nu})=\prod_{i=1}^{e_\nu} (1-\alpha_i T)\in
\zp[T]\] and put
\[\mathcal{P}_\nu=\mathcal{P}_\nu(A/k)=P_\nu (q_\nu^{-1}\gamma_\nu)\in
\La(\Gamma_\nu)\subseteq \La(\Gamma).\] If  one identifies
$\La(\Gamma_\nu)$ and $\La(\Gamma)$ with the power series rings
$\zp[[T_\nu]],$ $T_\nu=\gamma_\nu -1,$ and $\zp[[T]],$
$T=\gamma-1,$ for a fixed generator $\gamma$ of $\Gamma,$ then
$\mathcal{P}_\nu$ corresponds to
\[P_\nu(q_{\nu}^{-1}(T_\nu+1))=P_\nu(q_{\nu}^{-1}(T+1)^{f_\nu}),\]
where $f_\nu\in \zp$ is uniquely determined by the condition
\[q_\nu \omega(q_\nu^{-1})=\chi_{cyc}(\gamma)^{f_\nu}.\] Here
$\omega:\zp^\times\to\mu_{p-1}\subseteq \zp^\times$ and
$\chi_{cyc}:\Gamma\to 1+p\zp\subseteq \zp^\times$ denote the
Teichm\"{u}ller and cyclotomic character, respectively.

For a finitely generated $\La(G)$-module $M$ let
$\mathrm{char}_G(M)$ denote its characteristic ideal in $\La(G),$
where we assume $G\cong\zp.$ From the above considerations one
obtains immediately the following

\begin{prop} {\rm (}cf.\ \cite[prop.\ 2.4]{greenvat}{\rm )}\label{prop-euler}
\begin{enumerate}
\item
$\mathrm{char}_{\Gamma_\nu}(\H^1(k_{cyc,v},A)^\vee)=\La(\Gamma_\nu)\mathcal{P}_\nu$,
\item $\mathrm{char}_\Gamma (\Ind_\Gamma^{\Gamma_\nu}
\H^1(k_{cyc,v},A)^\vee)=\La(\Gamma)\mathcal{P}_\nu,$
\item The $\mu$-invariants of $
\H^1(k_{cyc,\nu},A)^\vee$ and $\Ind_\Gamma^{\Gamma_\nu}
\H^1(k_{cyc,\nu},A)^\vee$ are zero.
\item The $\lambda$-invariant $\lambda(\Ind_\Gamma^{\Gamma_\nu}
\H^1(k_{cyc,v},A)^\vee)$ is equal to $s_\nu d_\nu,$ where
$s_\nu=(\Gamma:\Gamma_\nu)$  equals $[(k_{cyc} \cap
k(\mu_{t})):k]$ with $t$  the largest power of $p$ dividing
$(q_\nu^{p-1}-1)$ and $ d_\nu=\lambda( \H^1(k_{cyc,v},A)^\vee)$ is
the multiplicity of $1-\widetilde{q_\nu^{-1}}T$ in
$\widetilde{P_\nu(T)}\in \mathbb{F}_p [T].$ Here $\widetilde{\cdot}$
means reduction modulo $p.$
\end{enumerate}
\end{prop}

\Subsection{The characteristic element of an elliptic curve over $k_\infty$}
Now we come to our arithmetic main results. Assume that $E$ is an elliptic
curve over $k$ with good ordinary reduction at all places $S_p.$

{\it Throughout the whole paper we assume that
 $E$ has good reduction at all places in $S_p.$}

As usual the $p$-Selmer group of $E$   is  defined as

\begin{eqnarray*}
\Sel_{p^\infty}(E/L)&:=&\ker\Big(  H^1(L, E_{p^\infty})\to
\bigoplus_{w}H^1(L_w,E(\overline{L_w}))_{p^\infty}\Big)\\
&\cong &\ker\Big(  H^1(G_S(L), E_{p^\infty})\to \bigoplus_{w\in S(L)}
H^1(L_w,E(\overline{L_w}))_{p^\infty}\Big).
\end{eqnarray*}
Here, $L$ is a finite extension of $k$ and, in the first line, $w$
runs through all places of $L$ while, in the second line, $S(L)$
denotes the set of all places of $L$
 lying above some place of $S.$ As usual, $L_w$ denotes the completion
 of $L$ at the place $w$ and for any field $K$ we fix an algebraic
 closure $\bar{K}.$ For infinite extensions $K$ of $k,$
$\Sel_{p^\infty}(E/K)$ is defined to be the direct limit of
$\Sel_{p^\infty}(E/L)$ over all finite intermediate extensions
$L.$

Suppose now that $k_\infty|k$ is an torsionfree pro-$p$ $p$-adic
Lie extension inside $k_S$ and  containing $k_{cyc}.$ Then its
Galois group is isomorphic to the semidirect product
$G:=G(k_\infty/k)\cong H\rtimes \Gamma$ where $H$ denotes the
Galois group of $k_\infty|k_{cyc}$ and $\Gamma=G(k_{cyc}/k)$ as
before.

The Selmer group $\Sel_{p^\infty}(E/k_\infty)$ bears a natural
structure as an discrete (left) $G$-module. For some purposes it
is more convenient to deal with (left) compact $G$-modules, thus
we take the Pontryagin duals $-^\vee$ and set
\[X(k_\infty):=(\Sel_{p^\infty}(E/k_\infty))^\vee.\]

\Subsection{The false Tate curve case}\label{tate-curve}

We first consider the case where $G$ is $2$-dimensional, i.e.\
isomorphic to the semidirect product $\zp\rtimes\zp.$
In this case Theorem \ref{ore-thm} tells us that
 $\mathcal{T}$ is an Ore set.


For simplicity we shall assume that $k$ contains the $p^\mathrm{th}$ roots
of unity. Let \linebreak $\mathfrak{M}_0(k_\infty/k)$ be a set of all primes
of $k$ which are not lying above $p$ and are ramified for
$k_\infty/k_{cyc}$. We put
\begin{multline}\label{p1}
\mathfrak{M}_1(k_\infty/k,E):=
\{v\in \mathfrak{M}_0(k_\infty/k)|\text{ $E/k$ has }\\
\text{split multiplicative reduction at  $v$}\},
\end{multline}
\begin{multline}
\label{p2} \mathfrak{M}_2(k_\infty/k,E):=
\{v\in \mathfrak{M}_0(k_\infty/k)|\text{ $E/k$ has good reduction}\\
\text{ at  $v$  and $E(k_v)_{p^\infty}\ne 0$.}\}
\end{multline}
and $\mathfrak{M}=\mathfrak{M}(k_\infty/k,E)
:=\mathfrak{M}_1(k_\infty/k,E)\cup \mathfrak{M}_2(k_\infty/k,E)$.

%
%

The following result, which relies heavily on the vanishing of higher $H$-homology groups of $X(k_\infty),$ generalizes partly the Euler characteristic
formula \cite[thm.\ 4.11]{hachi-ven}. More precisely the mentioned
formula will be reobtained in the Corollary below by ``evaluating the characteristic power
series $\mathrm{Ak}_G(H,X(k_\infty))$ at $0$" and applying the
$p$-adic valuation.

\begin{thm}\label{tate-case}
Assume that $X(k_{cyc})$ is a torsion $\La(\Gamma)$-module with
vanishing $\mu$-invariant. Then $X(k_\infty)$ is finitely
generated over $\La(H)$ and it holds modulo $\La(\Gamma)^\times$
that\[\pi_H(\mathrm{char}_G(X(k_\infty)))
\equiv \mathrm{Ak}_G(H,X(k_\infty))\equiv
\mathrm{char}_\Gamma(X(k_{cyc})) \cdot \prod_{\nu\in
\mathfrak{M}} \mathcal{P}_\nu(E(p)/k), \] where the local factors
are those defined in section \ref{local-euler}. 
\end{thm}

 Before we evaluate at ``$0$" we have to introduce some more notation.
We define the $p$-Birch-Swinnerton-Dyer constant as
\begin{equation*}
\rho_p(E/k):=
\frac{\sharp \sha(E/k)_{p^\infty}}
{(\sharp E(k)_{p^\infty})^2
\prod_{v}|c_v|_p}\times
\prod_{v|p}(\sharp\tilde{E_v}(\kappa_v)_{p^\infty})^2.
\end{equation*}
Here, $\sha(E/k)$ is the Tate-Shafarevich group of $E$ over $k$,
$\kappa_v$ is the residue field of $k$ at $v$ and $\Tilde{E_v}$ is
the reduction of $E$ over $\kappa_v$. We denote by $c_v$ the local
Tamagawa factor at $v$, $[E(k_v):E_0(k_v)]$, where $E_0(k_v)$ is
the subgroup of $E(k_v)$ consisting from all of the points which
maps to smooth points by reduction modulo $v$. $|*|_p$ denotes the
$p$-adic valuation normalized such that $|p|_p=\frac{1}{p}.$ For
any prime $v$ of $k$, let $L_v(E,s)$ be the local L-factor of $E$
at $v$. Let $P_0(k_\infty/k)$ be the set of all primes of $k$
which are not lying above $p$ and ramified for $k_\infty/K_{cyc}$.
As mentioned above  using Proposition \ref{EC-evaluation}) we (re)obtain

\begin{cor}\label{BSD}
In the situation of the theorem and assuming that the
$G$-Euler characteristic $\chi(G,X(k_\infty))$ of
$X(k_\infty)$ is finite let $F_{X(k_\infty)}\in \La_\mathcal{T}$ be a
characteristic element of $X(k_\infty).$ Then $F_{X(k_\infty)}(0)$ is
defined and non-zero, and it holds that
\[\chi(G,X(k_\infty))=|F_{X(k_\infty)}(0)|_p^{-1}=\rho_p(E/k)\times\prod_{v\in
\mathfrak{M}} |L_v(E,1)|_p.\]
\end{cor}


Before we give the proof of the theorem we introduce the modified Selmer group 

\begin{equation*}
\Sel'_{p^\infty}(E/K_{cyc}):=\mathrm{Ker} (H^1(k_S/k_{cyc},
E_{p^\infty})\rightarrow
\bigoplus_{S\backslash\mathfrak{M}}J_\nu(k_{cyc})),
\end{equation*}
where $ J_\nu(k_{cyc}))$ is the Pontryagin dual of
$\Ind_\Gamma^{\Gamma_\nu} \H^1(k_{cyc,\nu},E(p))^\vee,$ see
section \ref{local-euler}. Then we have the following exact
sequence
\begin{equation}\label{res-seq}
0\rightarrow \Sel_{p^\infty}(E/k_{cyc})\rightarrow
\Sel'_{p^\infty}(E/k_{cyc})\rightarrow
\bigoplus_{\mathfrak{M}}J_\nu(k_{cyc}) \rightarrow0.
\end{equation}

\begin{proof}
First note that $\H^i(H,\Sel_{p^\infty}(E/k_\infty))=0$ for $i\geq
1$ by the proof of \cite[thm.\ 4.11]{hachi-ven}. Thus
$\mathrm{Ak}_G(H,X(k_\infty))=\mathrm{char}_\Gamma(X(k_\infty)_H).$
But since according to the proof of (loc.\ cit.) the dual of the
restriction map
\begin{equation*}
\text{res}:\Sel'_{p^\infty}(E/k_{cyc})\rightarrow
\Sel_{p^\infty}(E/k_\infty)^H
\end{equation*}
is  a pseudo-isomorphism, the statement follows from the short
exact sequence \ref{res-seq}, the determination of the local
factors in Proposition \ref{prop-euler} and by proposition \ref{descent}.
\end{proof}

\begin{rem}\label{truncated}
If the rank of the Mordell-Weil group is positive we expect the vanishing of $F_{X(k_\infty)}(0)$ by the Birch and Swinnerton-Dyer conjecture. In this case one has to modify the Euler characteristic as is well known in the cyclotomic situation (\cite{perrin93} or \cite{schneider83},\cite{schneider85}) and as was proposed in \cite{cs2} in the $GL_2$-case. For simplicity we assume that $M$ is a $\La(G)$-module with $\H_i(H,M)=0$ for all $i\geq 1.$ Let \[\phi_M:\H_1(G,M)\to\H_0(G,M)\] be the following composition of maps \[\xymatrix@1@+10pt{
   {\ \H_1(G,M) \ar[r]^-{inf}_-{\cong} } &  {\ \H_1(\Gamma,\H_0(H,M)=(M_H)^\Gamma \ar[r]^-{\psi_M} } &  {\ (M_H)_\Gamma=\H_0(G,M),  }  
}\]
where $\psi_M:(M_H)^\Gamma\to (M_H)_\Gamma$ is induced by the identity on $M_H.$ We say that $M$ has  {\em finite truncated $G$-Euler characteristic,} if both $\coker(\phi_M)$ and $\ker(\phi_M)$ are finite,  and we define the truncated $G$-Euler characteristic of $M$ by \[\chi_t(G,M)=\#\coker(\phi_M)/\#\ker(\phi_M).\] Setting formally $H=1,$ e.g.\ $G=\Gamma,$ in the above we reobtain the definition of the generalized  $\Gamma$-Euler characteristic $\chi_t(\Gamma,N)$ of a $\La(\Gamma)$-module $N.$ Recall from \cite[lem.\ 3]{schneider83} that if $\chi_t(\Gamma,N)$ is defined (``semi-simplicity at zero"), then it equals $|c(N)|_p^{-1}$ where $c(N):=c(f_N):=[f_N(t)\cdot t^{-m(N)}]_{t=0}\in\O$ denotes the leading coefficient of the characteristic power series $f_N$ of $N$ (if the multiplicity of zero of $f_N(t)$ at $0$ is $m(N)$). Similarly, we introduce the {\em leading coefficient} $c(F)$ for any $F\in\La_\T$ by \[c(F):=c(\pi_H(F)).\] Now under the assumptions of Theorem \ref{tate-case} we have the following: If $\chi_t(G,X(k_\infty))$ is defined then \[\chi_t(G,X(k_\infty)=|c(F_{X(k_\infty)})|_p=\chi_t(\Gamma,X(k_{cyc}))\times\prod_{v\in
\mathfrak{M}} |L_v(E,1)|_p.\] This follows immediately from \cite[thm.\ 4.10]{hachi-ven} and its proof.
\end{rem}

\begin{prop}
In the situation of the theorem let $F_{X(k_\infty)}\in \La_\mathcal{T}$ be a
characteristic element of $X(k_\infty).$ If $\rho_G\to \mathrm{Aut}(V)$ is a finite dimensional representation of $G$ over a finite extension $K$ of $\qp$ such that   the $G$-Euler characteristic $\chi(G,V,X(k_\infty))$ of
$X(k_\infty)$ is finite, then $F_{X(k_\infty)}(\rho)$ is
defined and non-zero, and it holds that
\[\chi(G,V,X(k_\infty))=|F_{X(k_\infty)}(\rho)|_p^{-[K:\qp].}\] 
\end{prop}

\Subsection{The $GL_2$-case}\label{GL2}

Now suppose that $k_\infty=k(E(p))$ and hence $G=H\rtimes\Gamma,$
where $H$ is an open subgroup of $SL_2(\zp)$ and thus of dimension
$3.$ Assume that $H$ is uniform and that the semidirect-product is in fact direct, this can always be achieved after a finite base change. Then by Theorem \ref{ore-thm} the multiplicative set $\T=\La\setminus \m(H)$ is an Ore set and we can apply the theory of characteristic elements.

\begin{thm}\label{GL2-case}
Assume that $X(k_{cyc})$ is a torsion $\La(\Gamma)$-module with
vanishing $\mu$-invariant. Then $X(k_\infty)$ is finitely
generated over $\La(H)$ and it holds
that\[\pi_H(\mathrm{char}_G(X(k_\infty)))
\equiv \mathrm{Ak}_G(H,X(k_\infty))\equiv
\mathrm{char}_\Gamma(X(k_{cyc})) \cdot \prod_{\nu\in
\mathfrak{M}} \mathcal{P}_\nu(E(p)/k), \] where $\mathfrak{M}$ is
the set of places of $k$ with non-integral $j$-invariant and the
local factors are those defined in section \ref{local-euler}.
\end{thm}


This result generalizes partly \cite[thm.\ 3.1]{cs2}. Its proof is
analogous to that of theorem \ref{tate-case}, using now
\cite[rem.\ 2.6, lem.\ 3.3, lem.\ 3.6]{cs2}.
We leave it to the reader to derive an analogue of Corollary \ref{BSD} from this theorem. Also an analogue of Remark \ref{truncated} holds in this situation, see \cite[prop.\ 2.9, thm.\ 3.1]{cs2}.\\
We expect similar results for other $p$-adic Lie extensions of dimension at least $2$ and containing the cyclotomic $\zp$-extension and for other $p$-adic Galois representations. %
\POP
\PUSH{maincon.tex}%
\Section{Towards a main conjecture}\label{mainconj}

Let $L$ denote either the complex numbers $\mathbb{C}$ or a fixed
algebraic closure $\overline{\qp}$ of $\qp.$ By an {\em Artin
representation $\rho:G\to\mathrm{Aut}_L(V)$ over $L$} we mean a
finite dimensional representation of $G$ over $L$ which factorizes
through a finite quotient, say $\Delta,$ of $G.$ We fix embeddings
of a fixed algebraic closure $\bar{\Q}$ of $\Q$ both into
$\mathbb{C}$ and $\overline{\qp}$. Since an Artin representation
over $\mathbb{C}$ is already defined over a finite extension
$K\subseteq \bar{\Q}$ of $\Q$ we can also interpret   it as a
finite-dimensional $\overline{\qp}$-representation, and
vice-versa. Using character theory one sees immediately that the
equivalence classes of (absolutely irreducible) Artin
representations over $\mathbb{C}$ and $\overline{\qp}$ are
naturally equivalent. Let $E$ be an elliptic curve over $k.$  By $L(E,\rho,1)$ we denote the Hasse-Weil $L$-series, twisted by an Artin representation $\rho.$

Suppose as before that $k_\infty|k$ is an torsionfree pro-$p$ $p$-adic
Lie extension inside $k_S$ and  containing $k_{cyc}.$ Then its
Galois group is isomorphic to the semidirect product
$G:=G(k_\infty/k)\cong H\rtimes \Gamma$ where $H$ denotes the
Galois group of $k_\infty|k_{cyc}$ and $\Gamma=G(k_{cyc}/k)$ as
before.

As we have seen in subsection \ref{Ore-sets}, any $M\in\La(G)\mbox{-mod}^H$ gives rise to
an element of
 $K_0(\Lambda(G),\La(G)_{\T})$ and any element of
 the latter group  be represented by some $g\in
 (\La(G)_{\T})^\times.$ Assume that $\T$ is an Ore-set (e.g.\ if $H$ is uniform pro-$p$ and $G=H\times \Gamma$). We
 are quite optimistic and make the following

\begin{conj}\label{ver1}
Assume that $\T$ is an Ore-set. Let $E$ be an
elliptic curve over $k$ with good ordinary reduction at $S_p$ and
assume that $X(k_{cyc})$ is a $\La(\Gamma)$-torsion module with
vanishing $\mu$-invariant.  Then \begin{enumerate}
\item {\rm(existence of distribution)} there exist $F\in\La(G)_\T^\times$ such that for all
irreducible Artin-representations $\rho$  of $G$  such that
$L(E,\rho,1)$  is defined, i.e.\ $L(E,\rho,s)$ has no pole at
$s=1,$  $F(\rho)$ is defined and
\[F(\rho)=C(\rho)\cdot \prod_\M\mathrm{Euler}_\nu(E,\rho,1) \cdot \frac{L(E,\rho,1)}{\Omega_E(\rho)},\]
where $C(\rho)$ should be a generalization of Gau{\ss}-sums times a contribution from the Euler factors of $L(E,\rho,1)$ at primes above $p,$ $\mathrm{Euler}_\nu(E,\rho,1)$ the local   Euler factors of $L(E,\rho,1)$ at primes in $\M$ and $\Omega_E(\rho)$ a (Deligne)-Period.
\item {\rm(main conjecture)} $F$ is a characteristic element of $X(k_\infty).$
 \end{enumerate}
\end{conj}

\begin{rem}
 If a distribution $F $ with
the interpolation property in (i) exists, then  $F^\iota
$ interpolates $C(\rho^d)\cdot
\prod_\M\mathrm{Euler}_\nu(E,\rho^d,1)  \cdot
\frac{L(E,\rho^d,1)}{\Omega_E(\rho^d)}.$
\end{rem}

Of course, this is nothing but a proposal for the shape of a main conjecture and it will only be complete once we have a good guess for the occurring (epsilon-)factors, periods (depending on $\rho$). For the cyclotomic $\zp$-extension and motives over $\Q$ this has been discussed by J. Coates \cite{coates89}. We hope that his formalism can be adapted to our situation and we will come back to this item in a subsequent paper. Also we are aware that at moment our formalism applies only under restrictions on the base field, in particular we cannot take $\Q$ as base field in all our examples we are interested in. In order to illustrate  a possible setting consider the following example, already discussed in \cite{coates-howsonII} and \cite[example 8.7]{css}.

\begin{ex}
Let $E$ be the elliptic curve $X_1(11),$ given by $E\; :\; y^2+y=x^3-x^2,$ of conductor $11.$ Take $p=5,$  put $k=\Q(\mu_5),$ $k_\infty=k(E_{5^\infty})$ and $G=G(k_\infty/k),$ $H=G(k_\infty/\Q(\mu_{5^\infty})$ as well as $\Gamma=G(\Q(\mu_{5^\infty})/k).$ It is shown in \cite{fisher} that $G$ is isomorphic to the congruence subgroup \[\Gamma(5):=\{g\in GL_2(\mathbb{Z}_5)| g \mbox{ is congruent the identity modulo } 5\}.\] We have $G\cong H\times\Gamma$ where we have identified the center of $G$ with the quotient $\Gamma$ and we point out that $H$ is uniform. Thus, in this example $\T$ is known to be an Ore-set. Since $G(\Q(E_{5^\infty})/\Q)/H$ has the simple form $\Gamma\times\Delta$ where $\Delta$ is a finite group of order prime to $5$ we expect that one can even localize the Iwasawa algebra $\La( G(\Q(E_{5^\infty})/\Q))$ of the non pro-$5$ group $G(\Q(E_{5^\infty})/\Q)$ suitably. Then it would be possible to formulate a main conjecture over the base field $\Q$ as one  would  actually like. 
\end{ex}


A slightly stronger reformulation of the above conjecture reads as follows.
\newpage
\begin{conj}\label{vers2}
Let $E$ be an elliptic curve over $k$ with good ordinary reduction
at $S_p$ and assume that $X(k_{cyc})$ is a $\La(\Gamma)$-torsion
module with vanishing $\mu$-invariant.  Then \begin{enumerate}
\item {\rm (existence of distribution)} there exist $f,g\in\T$ such that for all
irreducible Artin-representations $\rho$  of $G$  such that
$L(E,\rho,1)$  is defined, i.e.\ $L(E,\rho,s)$ has no pole at
$s=1,$ one has $g(\rho)\neq 0$ and
\[f(\rho)g(\rho)^{-1}=C(\rho)\cdot \prod_\mathfrak{M}\mathrm{Euler}_\nu(E,\rho,1) \cdot \frac{L(E,\rho,1)}{\Omega_E(\rho)}.\]
\item {\rm(main conjecture)} we have the identity \[[X(k_\infty)]=[\Lambda/\Lambda f]-[\Lambda/\Lambda
g]\]
   in
$K_0(\Lambda(G)\mbox{-mod}^H).$
\end{enumerate}
\end{conj}

In view of Proposition \ref{chardet} and Theorem \ref{tate-case}
there should be a third equivalent version of the above conjectures
in the following style

\begin{conj}
Let $E$ be an elliptic curve over $k$ with good ordinary reduction
at $S_p$ and assume that $X(k_{cyc})$ is a $\La(\Gamma)$-torsion
module with vanishing $\mu$-invariant.  Then \begin{enumerate}
\item {\rm(existence of distribution)} there exist $F\in\La(G)_\T^\times$ such that for all
irreducible Artin-representations $\rho$  of $G$    one has
\[\mathrm{det}_{Q_\mathcal{O}(\Gamma)} (\pi_H (\mathrm{tw}_\rho(F))=\mathcal{L}(E\otimes\rho)\cdot \prod_{\nu\in
\mathfrak{M}} \mathcal{P}_\nu(E\otimes\rho/k), \] where
$\mathcal{L}(E\otimes\rho)$ denotes the (conjectural) cyclotomic
$p$-adic $L$-function associated with $E\otimes\rho$ and
$\mathcal{P}_\nu(E\otimes\rho/k)$ are the corresponding Euler
factors.
\item {\rm(main conjecture)} $F$ is a characteristic element of $X(k_\infty).$
\end{enumerate}
\end{conj}




A.\ Huber and G.\ Kings propose in \cite{hu-ki} an (non-commutative) Iwasawa main conjecture for motives from the point of view of the equivariant Tamagawa number (TNC) conjecture. Indeed, the validity of their main conjecture is equivalent to the validity of the equivariant TNC at each level of the tower of number fields. We are convinced that our Conjecture is coherent with their main conjecture (in cases where both conjectures are ``defined", also Huber-Kings consider the TNC specialized to the motive associated with $E$ only ``away from the critical point" and thus one has to consider an analogous version of their conjecture at this point). After a first version of this article was finished, T.\ Fukaya and K.\ Kato \cite{fukaya-kato} formulated a  slightly different main conjecture as Huber-Kings, though in a similar spirit. Their approach could be considered as an intermediate version between that of Huber-Kings and ours because   on the one hand they derive their main conjecture also from Tamagawa number conjectures, using conjectural $\epsilon$-elements with non-commutative coefficient rings, and on the other hand they also construct a localized $K_1.$ But instead of localizing the ring $\La$ they localize in some sense a certain category of bounded complexes of $\La$-modules, a construction which works in full generality but which is less explicit. The comparison of both approaches with ours could   also lead to the determination of the precise  constants $C(\rho)$ etc.\ in the above formulation. 
\POP




\INPUT{../bib/xbib.bib}
\INPUT{habil.bbl} 

\bibliographystyle{amsplain}
\bibliography{../bib/xbib}

\def\Dbar{\leavevmode\lower.6ex\hbox to 0pt{\hskip-.23ex \accent"16\hss}D}
  \def\cfac#1{\ifmmode\setbox7\hbox{$\accent"5E#1$}\else
  \setbox7\hbox{\accent"5E#1}\penalty 10000\relax\fi\raise 1\ht7
  \hbox{\lower1.15ex\hbox to 1\wd7{\hss\accent"13\hss}}\penalty 10000
  \hskip-1\wd7\penalty 10000\box7}
  \def\cftil#1{\ifmmode\setbox7\hbox{$\accent"5E#1$}\else
  \setbox7\hbox{\accent"5E#1}\penalty 10000\relax\fi\raise 1\ht7
  \hbox{\lower1.15ex\hbox to 1\wd7{\hss\accent"7E\hss}}\penalty 10000
  \hskip-1\wd7\penalty 10000\box7} \def\Dbar{\leavevmode\lower.6ex\hbox to
  0pt{\hskip-.23ex \accent"16\hss}D}
  \def\cfac#1{\ifmmode\setbox7\hbox{$\accent"5E#1$}\else
  \setbox7\hbox{\accent"5E#1}\penalty 10000\relax\fi\raise 1\ht7
  \hbox{\lower1.15ex\hbox to 1\wd7{\hss\accent"13\hss}}\penalty 10000
  \hskip-1\wd7\penalty 10000\box7}
  \def\cftil#1{\ifmmode\setbox7\hbox{$\accent"5E#1$}\else
  \setbox7\hbox{\accent"5E#1}\penalty 10000\relax\fi\raise 1\ht7
  \hbox{\lower1.15ex\hbox to 1\wd7{\hss\accent"7E\hss}}\penalty 10000
  \hskip-1\wd7\penalty 10000\box7}
\providecommand{\bysame}{\leavevmode\hbox to3em{\hrulefill}\thinspace}
\providecommand{\MR}{\relax\ifhmode\unskip\space\fi MR }
\providecommand{\MRhref}[2]{%
  \href{http://www.ams.org/mathscinet-getitem?mr=#1}{#2}
}
\providecommand{\href}[2]{#2}
\begin{thebibliography}{10}

\bibitem{bass}
H.~Bass, \emph{Algebraic {$K$}-theory}, W. A. Benjamin, Inc., New
  York-Amsterdam, 1968.

\bibitem{bourbaki}
N.~Bourbaki, \emph{\'{E}l\'ements de math\'ematique. {F}asc. {X}{X}{X}{I}.
  {A}lg\`ebre commutative. {C}hapitre 7: {D}iviseurs}, Hermann, Paris, 1965.

\bibitem{brumer}
A.~Brumer, \emph{Pseudocompact algebras, profinite groups and class
  formations}, J. of Algebra \textbf{4} (1966), 442--470.

\bibitem{bf}
D.~Burns and M.~Flach, \emph{Tamagawa numbers for motives with
  (non-commutative) coefficients}, Doc. Math. \textbf{6} (2001), 501--570
  (electronic).

\bibitem{coates89}
J.~Coates, \emph{On $p$-adic $l$-functions}, Ast{\'e}risque \textbf{177-178}
  (1989), no.~Exp.\ No.\ 701, 33--59.

\bibitem{coates91}
J.~Coates, \emph{Motivic {$p$}-adic {$L$}-functions}, $L$-functions and
  arithmetic (Durham, 1989), London Math. Soc. Lecture Note Ser., vol. 153,
  Cambridge Univ. Press, Cambridge, 1991, pp.~141--172.

\bibitem{coates-howsonII}
J.~Coates and S.~Howson, \emph{{Euler characteristics and elliptic curves II}},
  J. Math. Soc. Japan \textbf{53} (2001), 175--235.

\bibitem{cs2}
J.~Coates, P.~Schneider, and R.~Sujatha, \emph{{Links between cyclotomic and
  $GL_2$ Iwasawa theory}}, to appear in Doc. Math.

\bibitem{css}
\bysame, \emph{Modules over {I}wasawa algebras}, J. Inst. Math. Jussieu
  \textbf{2} (2003), no.~1, 73--108.

\bibitem{cr}
C.~W. Curtis and I.~Reiner, \emph{Methods of representation theory. {V}ol.
  {I}}, John Wiley \& Sons Inc., New York, 1981.

\bibitem{burns}
{D. Burns}, \emph{{Iwasawa theory and $p$-adic Hodge theory over
  non-commutative algebras I}}, preprint (1999).

\bibitem{d}
P.~Deligne, \emph{Le d\'eterminant de la cohomologie}, Current trends in
  arithmetical algebraic geometry (Arcata, Calif., 1985), Contemp. Math.,
  vol.~67, Amer. Math. Soc., Providence, RI, 1987, pp.~93--177.

\bibitem{dsms2}
J.~D. Dixon, M.~P.~F. du~Sautoy, A.~Mann, and D.~Segal, \emph{Analytic
  pro-{$p$} groups}, 2nd ed., Cambridge Studies in Advanced Mathematics,
  vol.~61, Cambridge University Press, Cambridge, 1999.

\bibitem{fisher}
T.~Fisher, \emph{Descent calculations for the elliptic curves of conductor 11},
  Proc. London Math. Soc. (3) \textbf{86} (2003), no.~3, 583--606.

\bibitem{fukaya-kato}
T.~Fukaya and K.~Kato, \emph{{A formulation of conjectures on $p$-adic zeta
  functions in non-commutative Iwasawa theory}}, preprint (2003).

\bibitem{green}
R.~Greenberg, \emph{{On the structure of certain Galois groups}}, Invent. Math.
  \textbf{47} (1978), 85--99.

\bibitem{green-adic}
\bysame, \emph{{Iwasawa Theory for $p$-adic Representations}}, Advanced Studies
  in Pure Mathematics \textbf{17} (1989), 97--137.

\bibitem{greenvat}
R.~Greenberg and V.~Vatsal, \emph{On the {I}wasawa invariants of elliptic
  curves}, Invent. Math. \textbf{142} (2000), no.~1, 17--63.

\bibitem{hachi-ven}
Y.~Hachimori and O.~Venjakob, \emph{Completely faithful selmer groups over
  kummer extensions}, to appear in Doc. Math.

\bibitem{howson}
S.~Howson, \emph{{Iwasawa theory of Elliptic Curves for $p$-adic Lie
  extensions}}, Ph.D. thesis, University of Cambridge, July 1998.

\bibitem{howson2000}
\bysame, \emph{Euler characteristics as invariants of {I}wasawa modules}, Proc.
  London Math. Soc. (3) \textbf{85} (2002), no.~3, 634--658.

\bibitem{howson2001}
\bysame, \emph{Structure of central torsion {I}wasawa modules}, Bull. Soc.
  Math. France \textbf{130} (2002), no.~4, 507--535.

\bibitem{hu-ki}
A.~Huber and G.~Kings, \emph{Equivariant {B}loch-{K}ato conjecture and
  non-abelian {I}wasawa main conjecture}, Proceedings of the International
  Congress of Mathematicians, Vol. II (Beijing, 2002) (Beijing), Higher Ed.
  Press, 2002, pp.~149--162.

\bibitem{kato}
K.~Kato, \emph{{ Hodge Theory and Values of Zeta Functions of Modular Forms}},
  erscheint in: Ast{\'e}risque.

\bibitem{km}
F.~F. Knudsen and D.~Mumford, \emph{The projectivity of the moduli space of
  stable curves. {I}. {P}reliminaries on ``det'' and ``{D}iv''}, Math. Scand.
  \textbf{39} (1976), no.~1, 19--55.

\bibitem{koch}
H.~Koch, \emph{Galoissche {T}heorie der {$p$}-{E}rweiterungen},
  Springer-Verlag, Berlin, 1970.

\bibitem{la}
M.~Lazard, \emph{Groupes analytiques $p$-adiques}, Publ. Math. I.H.E.S.
  \textbf{26} (1965), 389--603.

\bibitem{li-o}
H.~Li and F.~van Oystaeyen, \emph{Zariskian filtrations}, Kluwer Academic
  Publishers, Dordrecht, 1996.

\bibitem{mc-rob}
J.~C. McConnell and J.~C. Robson, \emph{Noncommutative {N}oetherian rings},
  John Wiley \& Sons Ltd., Chichester, 1987, With the cooperation of L. W.
  Small, A Wiley-Interscience Publication.

\bibitem{nsw}
J.~Neukirch, A.~Schmidt, and K.~Wingberg, \emph{Cohomology of number fields},
  Grundlehren der mathematischen Wissenschaften, vol. 323, Springer, 2000.

\bibitem{Ne}
A.~Neumann, \emph{Completed group algebras without zero divisors}, Arch. Math.
  \textbf{51} (1988), 496--499.

\bibitem{ochi-ven}
Y.~Ochi and O.~Venjakob, \emph{On the structure of {S}elmer groups over
  {$p$}-adic {L}ie extensions}, J. Algebraic Geom. \textbf{11} (2002), no.~3,
  547--580.

\bibitem{ochi-ven2}
\bysame, \emph{{On the ranks of Iwasawa modules over $p$-adic Lie extensions}},
  Math. Proc. Cambridge Philos. Soc. \textbf{135} (2003), 25--43.

\bibitem{pass}
D.~S. Passman, \emph{The algebraic structure of group rings},
  Wiley-Interscience [John Wiley \& Sons], New York, 1977.

\bibitem{perrin93}
B.~Perrin-Riou, \emph{Th\'eorie d'{I}wasawa et hauteurs {$p$}-adiques (cas des
  vari\'et\'es ab\'eliennes)}, S\'eminaire de Th\'eorie des Nombres, Paris,
  1990--91, Progr. Math., vol. 108, Birkh\"auser Boston, Boston, MA, 1993,
  pp.~203--220.

\bibitem{rubin}
K.~Rubin, \emph{{On the main conjecture of Iwasawa theory for imaginary
  quadratic fields.}}, Invent. Math. \textbf{93} (1988), no.~3, 701--713.

\bibitem{schneider83}
P.~Schneider, \emph{Iwasawa {$L$}-functions of varieties over algebraic number
  fields. {A} first approach}, Invent. Math. \textbf{71} (1983), no.~2,
  251--293.

\bibitem{schneider85}
\bysame, \emph{{p-adic height pairings. II.}}, Invent. Math. \textbf{79}
  (1985), 329--374.

\bibitem{sch-teit2}
P.~Schneider and J.~Teitelbaum, \emph{Algebras of $p$-adic distributions and
  admissible representations}, Preprintreihe des SFB 478 - Geometrische
  Strukturen in der Mathematik, M{\"u}nster \textbf{214} (2002).

\bibitem{serre65}
J.-P. Serre, \emph{Sur la dimension cohomologique des groupes profinis},
  Topology \textbf{3} (1965), 413--420.

\bibitem{smith}
P.~F. Smith, \emph{The {A}rtin-{R}ees property}, Paul Dubreil and Marie-Paule
  Malliavin Algebra Seminar, 34th Year (Paris, 1981), Lecture Notes in Math.,
  vol. 924, Springer, Berlin, 1982, pp.~197--240.

\bibitem{srinivas}
V.~Srinivas, \emph{Algebraic {$K$}-theory}, Progress in Mathematics, vol.~90,
  Birkh\"auser Boston Inc., Boston, MA, 1991.

\bibitem{swan}
R.~G. Swan, \emph{{Algebraic K-theory}}, LNM, vol.~76, Springer, 1968.

\bibitem{ven-weier}
O.~Venjakob, \emph{{A non-commutative Weierstrass preparation theorem and
  applications to Iwasawa theory}}, {J. reine angew. Math.} \textbf{{559}}
  ({2003}), 153--191.

\bibitem{ven}
\bysame, \emph{{Iwasawa Theory of $p$-adic Lie Extensions}}, Compos. Math.
  \textbf{138} (2003), no.~1, 1--54.

\bibitem{vo}
D.~Vogel, \emph{{Nonprincipal reflexive left ideals in Iwasawa algebras II}},
  preprint (2003).

\end{thebibliography}
\end{document}